\documentclass[10pt, a4paper]{article}

\usepackage{amsfonts, amsmath, amssymb, amsthm}
\usepackage{mathtools}
\usepackage{verbatim}
\usepackage{bbm}

\theoremstyle{plain}
\newtheorem{theorem}{Theorem}

\newtheorem{lemma}{Lemma}

\newtheorem{assumption}{Assumption}
\newtheorem{remark}{Remark}

\newtheorem{definition}{Definition}

\usepackage{algorithm}
\usepackage{algorithmic}
\usepackage{booktabs}
\usepackage{multirow}
\usepackage{longtable}
\usepackage{tabularx, makecell}

\usepackage[dvipsnames]{xcolor}
\usepackage{graphicx}

\allowdisplaybreaks

\usepackage{epstopdf}
\numberwithin{equation}{section}
\newcommand{\gbs}{\boldsymbol{g}}
\newcommand{\xbs}{{\boldsymbol{x}}}
\newcommand{\ybs}{{\boldsymbol{y}}}
\newcommand{\vbs}{{\boldsymbol{v}}}
\newcommand{\wbs}{{\boldsymbol{w}}}
\newcommand{\zbs}{{\boldsymbol{z}}}
\newcommand{\xref}{\xbs^\dagger}

\providecommand{\CB}{{\cal B}}
\newcommand{\CF}{{\cal F}}
\newcommand{\CD}{{\cal D}}
\newcommand{\bbE}{\mathbb{E}}
\newcommand{\bbN}{\mathbb{N}}
\newcommand{\bbR}{\mathbb{R}}
\newcommand{\bbP}{\mathbb{P}}
\newcommand{\goodset}{\mathcal{S}}
\newcommand*{\DP}[2]{\left<{#1},{#2}\right>}
\newcommand{\Null}{\operatorname{null}}						
\newcommand{\sign}{\mathrm{sign}}
\newcommand{\CJ}{{\cal J}}

\newcommand{\CX}{{\cal X}}
\newcommand{\CY}{{\cal Y}}
\newcommand{\VA}{{\mathbf{A}}}
\newcommand{\VB}{{\mathbf{B}}}
\newcommand{\dmapX}[1]{\CJ_{#1}^{\CX}}
\newcommand{\dmapXs}[1]{\CJ_{#1^\ast}^{\CX^{\ast}}}

\newcommand{\svaldmap}[1]{\jmath_{#1}}
\newcommand{\bregman}[2]{\VB_p(#1, #2)}
\newcommand{\bregmanS}[2]{\VB_{p^\star}\Big(\dmapX{p}(#1), \dmapX{p}(#2)\Big)}
\newcommand{\bregmanSdef}[2]{\VB_{p^\ast}\big(#1, #2\big)}
\newcommand{\norm}[1]{\|#1\|}
\newcommand{\xN}[1]{\|#1\|_{\CX}}
\newcommand{\xsN}[1]{\|#1\|_{\CX^\ast}}
\newcommand{\yN}[1]{\|#1\|_{\CY}}
\newcommand{\yNi}[1]{\|#1\|_{\CY_i}}

\newcommand{\ysNi}[1]{\|#1\|_{\CY^\ast_i}}

\newcommand{\Lmax}{L_{\max}}

\begin{document}
\title{Stochastic Gradient Descent for Nonlinear Inverse Problems in Banach Spaces\thanks{The work of B. Jin is supported by Hong Kong RGC General Research Fund (Projects 14306423 and 14306824),
ANR / RGC Joint Research Scheme (A-CUHK402/24) and a start-up fund from The Chinese University of Hong
Kong. Z. Kereta acknowledges support from the EPSRC (EP/X010740/1). The work of Y. Xia is supported by the National Natural Science Foundation of China (No. 12401556).}}
\date{}
\author{Bangti Jin\thanks{Department of Mathematics, The Chinese University of Hong Kong, Shatin, N.T., Hong Kong, China (\texttt{b.jin@cuhk.edu.hk})}
\and Zeljko Kereta\thanks{Computer Science Department, University College London, 169 Euston Road, London
NW1 2A, United Kingdom (\texttt{z.kereta@ucl.ac.uk})}
\and Yuxin Xia\thanks{Department of Mathematics, The Chinese University of Hong Kong, Shatin, N.T., Hong Kong, China (\texttt{yuxinxia@cuhk.edu.hk}); School of  Mathematics, Harbin Institute of Technology, Harbin 150001, Heilongjiang Province, China (\texttt{yuxinxia@hit.edu.cn})} }

\maketitle

\begin{abstract}
Stochastic gradient descent (SGD) and its variants are widely used and highly effective optimization methods in machine learning, especially for neural network training.
By using a single datum or a small subset of the data, selected randomly at each iteration, SGD scales well to problem size and has been shown to be effective for solving large-scale inverse problems.
In this work, we investigate SGD for solving nonlinear inverse problems in Banach spaces through the lens of iterative regularization.
Under general assumptions, we prove almost sure convergence of the iterates to the minimum distance solution and show the regularizing property in expectation under an \textit{a priori} stopping rule.
Further, we establish convergence rates under the conditional stability assumptions for both exact and noisy data.
Numerical experiments on Schlieren tomography and electrical impedance tomography are presented to show distinct features of the method.
\end{abstract}

\section{Introduction}\label{sec:intro}

This work is concerned with iterative methods for solving nonlinear inverse problems in
Banach spaces $\CX$ and $\CY$ that can be factorized into a system of $N$ nonlinear equations
\begin{equation}\label{eqn:inv_F_i}
F_i(\xbs)=\ybs_i, \quad i=1,\ldots, N,
\end{equation}
where each $F_i:\CD(F_i)\subset\CX\rightarrow\CY_i$ is a nonlinear forward operator between Banach spaces $\CX$ and $\CY_i$. We denote
$F=(F_1,\ldots, F_N)$, $\CY=\CY_1\otimes\cdots\otimes\CY_N$, and $\ybs=(\ybs_1,\ldots,\ybs_N)$.
In practice, due to the inevitable presence of measurement errors, we have access only to noisy data $\ybs^\delta = (\ybs_1^\delta, \ldots, \ybs_N^\delta)$,
which satisfies $\yNi{\ybs_i-\ybs_i^\delta}\le\delta_i$.
Let $\delta = \max\{\delta_1,\ldots,\delta_N\}$. Due to the ill-posedness of inverse problems, specialized techniques, such as regularization, have to be adopted.

Iterative regularization represents a powerful algorithmic paradigm and has been widely applied in many inverse problems. Classical iterative regularization methods include the Landweber method, Levenberg-Marquardt method, (iteratively regularized) Gauss-Newton method, and conjugate gradient method etc; see the monograph \cite{KaltenbacherNeubauerScherzer:2008} for an overview of regularizing properties of various iterative regularization methods in Hilbert spaces.
The main computational cost of iterative regularization methods lies in computing the gradient of the objective at each iteration, which can be prohibitively expensive for large-scale inverse problems.
Stochastic gradient descent (SGD) is a promising strategy to overcome this computational challenge.
SGD can be traced back to the seminal work of Robbins and Monro \cite{RobbinsMonro:1951}. The key idea is to use only a subset of the observational data $\ybs^\delta$ at each iteration, thereby reducing the computational cost per iteration.

SGD has had a profound impact on machine learning, particularly in the training of deep neural networks. In recent years, SGD and its variants have attracted growing interest within the inverse problems community and have been widely adopted in various imaging tasks, e.g., computed tomography \cite{HermanMeyer:1993,KarimiWard:2017}, positron emission tomography \cite{EhrhardtSchonlieb:2019,KeretaTwyman:2021,TwymanThielemans:2023} and optical tomography \cite{ChenLiLiu:2018,MacDonaldArridge:2020}.
The first work to analyze SGD through the lens of regularization theory is due to Jin and Lu \cite{JL19} in which the regularizing property of SGD was shown for linear inverse problems in Euclidean spaces.
Since then, substantial advancements have been made in both theoretical and practical aspects of SGD for inverse problems, including extensions for nonlinear inverse problems \cite{GFHF24,HuangJinLuZHang25,JinZhouZou:2020}, early stopping rules \cite{HuangJinLuZHang25,JahnJin:2020}, developments regarding step-size schedules \cite{JinLuZHang:2023,RabeloLeitao:2022}, data-driven approaches \cite{ZhouDSGD:2024} and Banach space formulations \cite{GFHF24,HuangJinLuZHang25,JinKereta:2023}; see the recent surveys \cite{EhrhardtKereta:2025, JinXiaZhou:2025} for further details.

Regularization methods in Banach spaces provide a flexible and principled framework to capture distinct features of the solutions.
In contrast to Euclidean and Hilbert spaces, which tend to favor smoothness, Banach spaces offer richer geometric structures.
This flexibility is essential for problems that are more adequately described in non-Hilbert settings. The sought solutions may have \textit{a priori} properties, e.g., sparsity and piecewise constancy, and the observed data may contain non-Gaussian noise.
These require changing either the solution space $\CX$, the data space $\CY$, or both.
For example, inverse problems with impulsive noise are more suitably modeled by choosing the data space $\CY$ as a Lebesgue space $L^r$ with $r\approx 1$ \cite{ClasonJinKunisch:2010}, while the recovery of sparse solutions is modeled by doing the same to the solution space $\CX$ \cite{CandesRombertTao:2006}.

The study of SGD in Banach spaces in the context of regularization theory is still in its infancy. Jin and Kereta \cite{JinKereta:2023} studied SGD for linear inverse problems in Banach spaces, and proved its regularizing property. In parallel, Jin et al. \cite{JinLuZHang:2023} proposed stochastic mirror descent (SMD) for linear inverse problems between a Banach space $\CX$ and a Hilbert space $\CY$, which can be interpreted as a randomized gradient method applied to the dual problem. These methods have been extended for solving nonlinear inverse problems \cite{GFHF24,HuangJinLuZHang25}.
Huang et al. \cite{HuangJinLuZHang25} proposed an \textit{a posteriori} stopping rule.
However, the works \cite{GFHF24,HuangJinLuZHang25} focus primarily on the convergence analysis and do not provide convergence rates, and use step-sizes that depend on the current iterate.
To the best of the authors’ knowledge, establishing convergence rates even for the classical Landweber method with noisy data between two Banach spaces $\CX$ and $\CY$, remains an open problem.

In this work, we investigate a general SGD method for nonlinear inverse problems in Banach spaces, and derive convergence rates for both exact and noisy data under a conditional stability estimate on the forward operator $F$, cf. Theorems \ref{thm:cond_stability_convergence_qeqp} and \ref{thm:noisy_convergence_rate}, which represents an important step forward in the theory of iterative regularization for inverse problems.
In addition, we establish its convergence and regularizing properties.
Note that, in the nonlinear setting, the analysis requires techniques distinct from those used in linear case, as it is necessary to ensure that the iterates remain within a region where the forward maps $F_i$ are well behaved (cf. Assumption \ref{assn:F}).
We first consider the case of exact data and show that SGD iterates converge to the $\xbs_0$-minimum-distance solution $\xref$  (first almost surely and then in expectation) under suitable assumptions on the nonlinear operators $F_i$, the step-size, and geometric properties of the space $\CX$; see Theorems \ref{thm:cone_convergence_qeqp} and \ref{thm:L1_cone_convergence_qeqp}.
For noisy observations, in Theorem \ref{thm:regularisation_property} we show the regularizing property of SGD for properly chosen stopping indices. In addition, we perform numerical experiments on Schlieren tomography and electrical impedance tomography to illustrate distinct features of the method.

The rest of the paper is organized as follows. In Section \ref{sec:prelims}, we collect standard notions for the geometry of Banach spaces, e.g.,  duality map and Bregman distance. In Section \ref{sec:nonlinear}, we describe the method, state the assumptions on the forward model and spaces $\CX$ and $\CY$, and present the convergence analysis.
In Section \ref{sec:nonlin_rate}, we establish convergence rates under conditional stability for both exact and noisy data. In Section \ref{sec:numerics}, we present numerical experiments on two image reconstruction problems. In appendices \ref{app:L1_cone_convergence_qeqp} and \ref{app:cone_convergence_qeqp} we give several technical tools and lengthy proofs. Throughout, let $\CX$ and $\CY$ be two real Banach spaces, with norms denoted by $\xN{\cdot}$ and $\yN{\cdot}$, respectively.
$\CX^\ast$ and $\CY^\ast$ are their respective dual spaces, with norms denoted by $\|\cdot\|_{\CX^\ast}$ and $\|\cdot\|_{\CY^\ast}$, respectively.
For $\xbs\in\CX$ and $\xbs^\ast\in\CX^\ast$, we denote the corresponding duality pairing by
$\DP{\xbs^\ast}{\xbs}=\DP{\xbs^\ast}{\xbs}_{\CX^\ast\times\CX}=\xbs^\ast(\xbs)$.
For a continuous linear operator $\VA:\CX\rightarrow \CY$, we use $\|\VA\|_{\CX\to \CY}$ to indicate the operator norm (often with the subscript omitted).
The adjoint of $\VA$ is denoted as $\VA^\ast:\CY^\ast\rightarrow\CX^\ast$ and it is a continuous linear operator, with $\norm{\VA}_{\CX\to\CY}=\norm{\VA^\ast}_{\CY^\ast\to\CX^\ast}$.
The conjugate exponent of $p\in (1,\infty)$ is denoted by $p^\ast$, such that $1/p+1/p^\ast=1$ holds.
The Cauchy-Schwarz inequality $| \DP{\xbs^\ast}{\xbs}|\le \xsN{\xbs^\ast}\xN{\xbs}$ holds for any $\xbs\in\CX$ and $\xbs^\ast\in\CX^\ast$.
For reals $a$ and $b$ we write $a\wedge b=\min\{a,b\}$ and $a\vee b=\max\{a,b\}$. By $(\CF_k)_{k\in\bbN}$,
we denote the natural filtration, i.e., a growing sequence of $\sigma$-algebras  such that $\CF_k
\subset\CF_{k+1}\subset \CF$, for all $k\in\bbN$ and a $\sigma$-algebra $\CF$, and $\CF_k$ is generated by
random indices $i_j$, for $j\leq k$.
In the context of SGD, $k\in\bbN$ is the iteration number and
$\CF_k$ denotes the iteration history, i.e., information available at time $k$, and for a given
initialization $\xbs_0$, we can identify $\CF_k=\sigma(\xbs_1,\ldots,\xbs_k)$.
For a filtration $(\CF_k)_{k\in\bbN}$, we denote by $\bbE_k[\cdot]=\bbE[\cdot\mid\xbs_1,\ldots,\xbs_k]$ the conditional expectation with respect to $\CF_k$. The notation a.s. denotes almost sure events.

\section{Preliminaries on Banach spaces}\label{sec:prelims}
In this section, we recall essential concepts from Banach space theory, including the duality maps and Bregman distances. For a comprehensive overview of Banach space geometry, see the monographs \cite{C_09,TBHK_12}.

We first recall standard notions of smoothness and convexity of Banach spaces.

\begin{definition}\label{defn:smoothness_and_convexity}
Let $\CX$ be a Banach space and define the function $\delta_\CX:(0,2]\rightarrow\bbR$ by
\begin{align*}
    \delta_\CX(\tau) = \inf\Big\{1-\frac{1}{2}\xN{\zbs+\wbs} : \xN{\zbs}=\xN{\wbs}=1, \xN{\zbs-\wbs}\geq\tau \Big\}.
\end{align*}
$\CX$ is said to be $p$-convex, for $p>1$, if $\delta_\CX(\tau)\geq K_p\tau^p$ for some $K_p>0$ and all $\tau\in(0,2]$.
The function $\rho_\CX:(0,\infty)\rightarrow(0,\infty)$ defined by
\begin{align*}
    \rho_\CX(\tau) = \sup\Big\{\frac{\xN{\zbs+\tau\wbs}+\xN{\zbs-\tau\wbs}}{2}-1 : \xN{\zbs}=\xN{\wbs}=1\Big\}
\end{align*}
is called the \emph{modulus of smoothness} of $\CX$, and is a convex and continuous function such that $\frac{\rho_\CX(\tau)}{\tau}$ is a non-decreasing function with $\rho_\CX(\tau)\leq\tau$.
The space $\CX$ is said to be \emph{uniformly smooth} if $\lim_{\tau\searrow0}\frac{\rho_\CX(\tau)}{\tau}=0$, and is said to be $p$-smooth, for $p>1$, if $\rho_\CX(\tau)\leq K_p\tau^p$ for some $K_p>0$ and all $\tau\in(0,2]$.
\end{definition}

Important examples of Banach spaces include Lebesgue spaces $L^p(\Omega)$, Sobolev spaces $W^{n,p}(\Omega)$ (for an open bounded domain $\Omega$, $n\in\bbN$), and the sequence spaces $\ell^p(\bbR)$ of $p$-summable real sequences. It can be shown that these Banach spaces are all $p\wedge 2$-smooth and $p\vee2$-convex for all $1<p<\infty$.

\begin{definition}[Duality map]\label{defn:duality_map}
For any $p>1$, a \emph{duality map} $\dmapX{p}:\CX\rightarrow 2^{\CX^\ast}$ is defined by
\begin{align*}
\dmapX{p} (\xbs) = \{\xbs^\ast\in\CX^\ast : \DP{\xbs^\ast}{\xbs}=\xN{\xbs}\xsN{\xbs^\ast}, \text{ and } \xN{\xbs}^{p-1}=\xsN{\xbs^\ast} \},
\end{align*}
with the gauge function $t\mapsto t^{p-1}$. A single-valued selection of $\dmapX{p}$ is denoted by $\svaldmap{p}$.
\end{definition}

For every $\xbs\in\CX$, the set $\dmapX{p}(\xbs)$ is non-empty, convex, and weakly closed in $\CX^\ast$.
The parameter $p$ typically depends on  geometric properties of $\CX$.
The following properties of Banach spaces and duality map hold.
\begin{theorem}[{\cite[Theorems 2.29, 2.52 and 2.53]{TBHK_12}}]\label{thm:dmap_properties}
\begin{enumerate}
\item[{\rm(i)}] The duality map $\dmapX{p}$ is monotone:
\begin{equation}\label{eqn:dual_non-negative}
    \DP{\svaldmap{p}(\xbs)-\svaldmap{p}(\ybs)}{\xbs-\ybs}\ge 0, \quad \forall \xbs,\ybs\in\CX, \,\, \svaldmap{p}(\xbs)\in \dmapX{p}(\xbs), \,\,\svaldmap{p}(\ybs)\in \dmapX{p}(\ybs).
\end{equation}
 Additionally, if $\CX$ is strictly convex, then  $\DP{\svaldmap{p}(\xbs)-\svaldmap{p}(\ybs)}{\xbs-\ybs}= 0 \Leftrightarrow \xbs=\ybs$.
\item[{\rm(ii)}]\label{rem:Xp_Xstarpstar}
$\CX$ is $p$-smooth if and only if $\CX^\ast$ is $p^\ast$-convex.
$\CX$ is $p$-convex if and only if $\CX^\ast$ is $p^\ast$-smooth.
\item[{\rm(iii)}]\label{property:dual_inverse}
$\CX$ is smooth if and only if every duality map $\dmapX{p}$ is single valued. If $\CX$ is smooth and reflexive, then $\dmapX{p}$ is invertible and $\big(\dmapX{p}\big)^{-1}=\dmapXs{p}$.
If $\CX$ is uniformly smooth and uniformly convex, then $\dmapX{p}$ and $\dmapXs{p}$ are both uniformly continuous on bounded sets.
\end{enumerate}
\end{theorem}

Due to the geometry of Banach spaces, it is often more convenient and appropriate to use the Bregman distance instead of the standard Banach space norm  in the analysis. In the rest of the paper, the dependence of $\bregman{\zbs}{\wbs}$ on the space $\CX$ is omitted, since it is often clear from context.
\begin{definition}\label{defn:Bregman_distance}
For $1/p+1/p^\ast=1$ and a smooth Banach space $\CX$, the functional
\begin{align*} \bregman{\zbs}{\wbs}
&= \frac{1}{p^\ast}\xN{\zbs}^p+\frac{1}{p}\xN{\wbs}^p-\DP{\dmapX{p}(\zbs)}{\wbs},
\end{align*}
is called the \emph{Bregman distance} with respect to the gauge function $\zbs\mapsto \xN{\zbs}^p$.
\end{definition}

We now present some useful  properties of the Bregman distance that depend on the geometry of the underlying Banach space and its duality map.
\begin{theorem}[{\cite[Theorem 2.60, Lemmas 2.62 and 2.63]{TBHK_12}}]\label{thm:bregman_properties}
\begin{enumerate}
\item[{\rm(i)}] If $\CX$ is smooth and reflexive, then Bregman distances in $\CX$ and $\CX^\ast$ satisfy
\[\bregman{\zbs}{\wbs}=\bregmanS{\wbs}{\zbs}.\]
\item[{\rm(ii)}] Bregman distance satisfies the three-point identity
\begin{align}\label{eqn:3_point_id}
\bregman{\zbs}{\wbs}=\bregman{\zbs}{\vbs}+\bregman{\vbs}{\wbs}+\DP{\dmapX{p}(\vbs)-\dmapX{p}(\zbs)}{\wbs-\vbs}.
\end{align}
\item[{\rm(iii)}]\label{rem:bregman_pconv} If $\CX$ is $p$-convex, then it is reflexive, $p\geq 2$ and there exists $C_p>0$ such that
\begin{align}\label{eqn:norm_leq_bregman}
\bregman{\zbs}{\wbs}\geq\frac{C_p}{p}\xN{\wbs-\zbs}^p.
\end{align}
If $\CX^\ast$ is $p^\ast$-smooth, then it is reflexive, $p^\ast\le 2$ and there exists $G_{p^\ast}>0$ such that
\begin{align*}
\bregmanSdef{\zbs^\ast}{\wbs^\ast}\leq\frac{G_{p^\ast}}{p^\ast}\xsN{\wbs^\ast-\zbs^\ast}^{p^\ast}.
\end{align*}
\item[{\rm(iv)}] $\bregman{\zbs}{\wbs}$ is continuous in the second argument. If $\CX$ is smooth and uniformly convex, then $\dmapX{p}$ is continuous on bounded subsets and $\bregman{\zbs}{\wbs}$ is continuous in its first argument.
\end{enumerate}
\end{theorem}

\begin{lemma}[{Coercivity of the Bregman distance \cite[Lemma A.3]{JinKereta:2023}}]\label{lem:bregman_coercivity}
If $\bregman{\xbs}{\xref}\leq C<\infty$, then $\xN{\xbs}^p\leq (2p^\ast)^{p}(\xN{\xref}^p\vee C)$.
\end{lemma}

\section{Convergence analysis}\label{sec:nonlinear}

In this section, we analyze the convergence and regularizing properties of the SGD. Throughout we assume that there exists a solution $\xref$ such that $F(\xref)=\ybs$, which implies $\cap_{i=1}^N\CD(F_i)=\CD(F)\neq\emptyset$. We denote by $\CX_{\min}$ the set $\{\xbs\in \CX: F(\xbs)=\ybs\}$.

\begin{remark}\label{rmk:product}
The set $\CY=\otimes_{i=1}^N\CY_i$ is equipped with the $\ell^r$ norm, for $r\geq 1$, defined via
\begin{align}\label{eqn:lr_norm}
    \|{\ybs}\|_{\CY,\ell^r}:= \|(\|\ybs_1\|_{\CY_1},\ldots,\|\ybs_N\|_{\CY_N})\|_{r}.
\end{align}
For $r,q\in[1,\infty]$, the norms $\|{\ybs}\|_{\CY,\ell^r}$ and $\|{\ybs}\|_{\CY,\ell^q}$ are equivalent, i.e.,
$\|{\ybs}\|_{\CY,\ell^r}^q \leq C_{q,r}\|{\ybs}\|_{\CY,\ell^q}^q$, with $C_{q,r}=1$ if $r\ge q$ and $C_{q,r}=N^{q/r-1}$ if $q>r$.
Below we omit the subscript $r$, and denote the constant by $C_q$.
\end{remark}

To construct an approximate solution of problem \eqref{eqn:inv_F_i}, from possibly noisy measurements $\ybs^\delta$, consider the following objective
\begin{equation}\label{eqn:obj_function}
\Psi(\xbs)=\frac{1}{N}\sum_{i=1}^N \Psi_i(\xbs), \quad \mbox{with }\Psi_i(\xbs)=\frac{1}{q}\yNi{F_i(\xbs)-\ybs_i^\delta}^q,  \quad \xbs\in \CD(F).
\end{equation}
To get a minimizer of $\Psi(\xbs)$,  given the initial guess $\xbs_0$, consider the following iteration
\begin{align}\label{eqn:sgd}
\xbs_{k+1}^\delta = \dmapXs{p}\left(\dmapX{p}(\xbs_k^\delta) - \mu_{k+1}\gbs_{k+1}^\delta\right), \quad k=0,1,\ldots,
\end{align}
where $\mu_{k+1}>0$ is the step-size and $\gbs_{k+1}^\delta=\gbs(\xbs_k^\delta,\ybs^\delta,i_{k+1})$ is the stochastic update direction. For exact data $\ybs$, we will omit the superscript $\delta$ from relevant quantities.
The choice of the stochastic function $\gbs$ defines the given stochastic gradient method.
For SGD, we define
\begin{align*}
    \gbs(\xbs,\ybs^\delta,i) = F^\prime_i(\xbs)^\ast\svaldmap{q}(F_i(\xbs)-\ybs_i^\delta)=\partial\Psi_i(\xbs),
    \quad \mbox{for }1<q\leq p.
\end{align*}
The index $i_k$ is chosen uniformly at random from the set $\{1,\ldots, N\}$.
The stochastic gradient is an unbiased estimator of the full sub-gradient $\partial\Psi(\xbs)$, since $\bbE[\gbs(\xbs,\ybs^\delta,i)] = \partial\Psi(\xbs)$.
Using $\gbs(\xbs,\ybs^\delta,i_k)$, instead of $\partial\Psi(\xbs)$, reduces the per-iteration cost roughly by a factor of $N$.

The analysis requires suitable assumptions on Banach spaces $\CX$ and $\CY$ and the operators $F_i$, which we state next.
The assumption on $\CY$ is only required to prove stability results in Section \ref{sec:reg_property}.
\begin{assumption}\label{assn:X_and_Y}
$\CX$ is $p$-convex and uniformly smooth, and $\CY$ is uniformly smooth.
\end{assumption}

For a general nonlinear mapping $F$, the objective $\Psi$ is non-convex, which poses significant challenges in the analysis of gradient-based methods for the solution.
To ensure the (local) convergence of an iterative method, we require (local) properties of the forward operator $F$, which provide sufficient control of the degree of nonlinearity around the target solution $\xref$.
We denote a ball of radius $\nu>0$ in the Bregman distance by $\CB_\nu(\xbs)=\{\zbs\in\CX: \bregman{\zbs}{\xbs}\leq \nu\}$.

\begin{assumption}\label{assn:F}
Let $\goodset\subset \CD(F)$ be a closed and convex set.
The forward operators $F_i:\CD(F_i)\subset\CX\to \CY_i$, for $i=1,\dots, N$, satisfy the following conditions.
\begin{enumerate}
    \item[{\rm(i)}] There exists a solution $\xref\in\CX$ of \eqref{eqn:inv_F_i} and an $\nu>0$ such that $\CB_{\nu}(\xref)\subset \goodset \subset \CD(F)$.
    \item[{\rm(ii)}] The operators $F_i$ are Fr\'{e}chet differentiable in $\goodset$ and the their derivatives satisfy $\|F^\prime_i(\xbs)\|\leq L_i$, for all $\xbs\in\goodset $. Let $\Lmax:=\max_{i\in[N]} L_i$.
    \item[{\rm(iii)}] There exists $0<\gamma<1$ such that the local tangential cone condition holds
\begin{align}\label{eqn:TCC}
\yNi{F_i(\xbs)-F_i(\tilde\xbs)-F_i'(\xbs)(\xbs-\tilde\xbs)}\leq \gamma \yNi{F_i(\xbs)-F_i(\tilde\xbs)},\quad
\forall \xbs,\tilde\xbs\in\goodset.
\end{align}
\end{enumerate}
\end{assumption}

Assumption \ref{assn:F} is frequently employed in the analysis of iterative regularization methods for nonlinear inverse problems. 
The set $\goodset$ in (i) and (iii) does not play a role for the moment and can be chosen as $\CB_{\nu}(\xref)$.
However, in Lemma \ref{lem:unique_mini_dis_solu} below, additional conditions will be imposed on $\goodset$ to ensure the existence and uniqueness of an $\xbs_0$-minimum-distance solution.
The tangential cone condition in (iii) controls the degree of the nonlinearity of  $F_i$, and has been extensively used  \cite[Section 7.1.3]{TBHK_12}.

\begin{remark}\label{rem:TCC_frech_obj}
If the region of validity of \eqref{eqn:TCC} is convex, then \eqref{eqn:TCC} enforces the convexity of the set $\CX_{\min}$.
Moreover, the tangential cone condition \eqref{eqn:TCC} yields a useful inequality
\begin{align}\label{eqn:TCC_1}
(1-\gamma)\yNi{F_i(\xbs)-F_i(\tilde\xbs)}    \leq \yNi{F_i'(\xbs)(\xbs-\tilde\xbs)}\leq (1+\gamma) \yNi{F_i(\xbs)-F_i(\tilde\xbs)},\quad\forall
\xbs,\tilde\xbs\in\goodset.
\end{align}
The $L_i$-boundedness of the derivative $F^\prime_i$ of $F_i$ and \eqref{eqn:TCC_1} give a local Lipschitz bound on $F_i$ as
\begin{align*}
\yNi{F_i(\xbs)-F_i(\tilde\xbs)}\leq \tfrac{L_i}{1-\gamma}\xN{\xbs-\tilde\xbs}, \quad \forall \xbs,\tilde\xbs\in\CB_\nu(\xref).
\end{align*}
\end{remark}

\begin{lemma}
Let Assumption \ref{assn:F} hold and $\xbs\in\CB_\nu(\xref)$.Then $\gbs(\xbs,\ybs^\delta,i)$ satisfies the moment bound
\begin{align}\label{eqn:variance_bound}
    \bbE[\xsN{\gbs(\xbs,\ybs^\delta,i)}^{p^\ast}] \leq \bbE[L_i^{p^\ast} \yNi{F_i(\xbs)-\ybs_i^\delta}^{(q-1)p^\ast}].
\end{align}
\end{lemma}
\begin{proof} Under the assumptions of the lemma, we have
\begin{align}\label{eqn:q_var_basic_bound}
    \xsN{\gbs(\xbs,\ybs^\delta,i)}
    &= \xsN{F^\prime_i(\xbs)^\ast\svaldmap{q}(F_i(\xbs)-\ybs_i^\delta)}\nonumber\\
    &\leq L_i \ysNi{\svaldmap{q}(F_i(\xbs)-\ybs_i^\delta)}
    \leq L_i \yNi{F_i(\xbs)-\ybs_i^\delta}^{q-1}.
\end{align}
Taking the expectation completes the proof.
\end{proof}

The central aim of this work is to show the convergence, and establish convergence rates, of iterates \eqref{eqn:sgd} towards a solution of \eqref{eqn:inv_F_i}. This requires a reference solution $\xref$, for which we adopt the so-called minimum distance solution.
\begin{definition}\label{defn:min_distance_solu}
Let $\goodset$ be the set in which $F$ satisfies Assumption \ref{assn:F}.
An $\xref\in\goodset$ satisfying
$$F(\xref)=\ybs\quad \mbox{and}\quad \bregman{\xbs_0}{\xref} =\inf\{\bregman{\xbs_0}{\xbs}: F(\xbs)=\ybs \text{ and } \xbs\in\goodset\}$$
is called the $\xbs_0$-minimum-distance solution.
\end{definition}

The next lemma states the existence and uniqueness of an $\xbs_0$-minimum-distance solution.
\begin{lemma}[{\cite[Lemma 7]{M18}}]\label{lem:unique_mini_dis_solu}
    Let Assumptions \ref{assn:X_and_Y} and \ref{assn:F} hold. Assume that $\goodset$ is closed and convex, and that $F$ is weak-to-weak continuous. Then an $\xbs_0$-minimum-distance solution exists and is unique. Moreover, any $\xbs^\ast$ such that $F(\xbs^\ast)=\ybs$ satisfies
    \begin{equation*}
        \DP{\dmapX{p}(\xbs_0)-\dmapX{p}(\xref)}{\xbs^\ast-\xref}\leq0.
    \end{equation*}
\end{lemma}

\begin{remark}
Let $\goodset$ be the set in which $F$ satisfies Assumption \ref{assn:F}. The $\xbs_0$-minimum-norm solution, defined as $\xbs^\dagger\in\goodset$ such that
$$F(\xref)=\ybs \quad\mbox{and}\quad \xN{\xbs_0-\xref} =\inf\{\xN{\xbs_0-\xbs}: F(\xbs)=\ybs\mbox{ and }\xbs\in\goodset\},$$
is also commonly employed \cite{JinKereta:2023,JinZhouZou:2020}. Its existence and uniqueness are shown in \cite[Proposition 7.1]{TBHK_12}. By suitably adapting the analysis, we can also show convergence results for an $\xbs_0$-minimum-norm solution; see Theorem \ref{thm:app_cone_convergence_qeqp} in the appendix.
\end{remark}

Now we briefly discuss the choice of the exponent $q$ in the objective $\Psi(\xbs)$. There are two notable special cases of $q$.
First, if $q=p$, because $(p-1)p^\ast=p$, the variance bound \eqref{eqn:variance_bound} can be stated in terms of the objective, i.e., $\bbE[\xsN{\gbs(\xbs,\ybs^\delta,i)}^{p^\ast}] \leq p\Lmax^{p^\ast} \Psi(\xbs)$, which simplifies the convergence analysis.
Second, if  $\CY_i=\ell^r$ and $q=r$, then from \eqref{eqn:lr_norm}, we have
$\frac{1}{qN}\sum_{i=1}^N\yNi{F_i(\xbs)-\ybs_i^\delta}^{q} = \frac{1}{qN}\yN{F(\xbs)-\ybs^\delta}^q$, which resembles the relationship between the full forward operator $F$ and sub-operators $F_i$ in Euclidean spaces when using the squared norm.
For noisy iterates, it is common to impose additional smoothness assumptions on $\CY$, since the duality map is generally continuous only at the origin. In such a case, for an $r$-smooth $\CY$, it might be desirable to pick $q=r$ instead of $p$. In this work, however, we focus on the case $q=p$.

The challenge in analyzing the convergence of the iterations \eqref{eqn:sgd} lies  in ensuring that the iterates stay within a region $\goodset$ in which Assumption \ref{assn:F} holds.
This is critical when dealing with nonlinear forward operators whose smoothness properties hold only in a local neighborhood, as the iterates $(\xbs_k)_{k\in\bbN}$ do not necessarily stay within a certain ball around the initial point.
Thus, the analysis differs markedly from the linear case.
In the case of exact data, proving that iterates $(\xbs_k)_{k\in\bbN}$ remain within a small ball  can be achieved by following the line of reasoning in the linear case \cite[Lemma A.2]{JinKereta:2023}, and using the tangential cone condition \eqref{eqn:TCC} to deal with the non-linearity of the forward operator.
However, for noisy data, such a result is no longer easily attainable unless specific step-size rules are adopted \cite{GFHF24,HuangJinLuZHang25}.

\subsection{Convergence analysis for exact data}\label{sec:convergence_property}
We now analyze the convergence of SGD for exact data (i.e., for $\delta=0$).
Note that under the tangential cone condition the following weak local convexity of the functional $\Psi(\xbs)$ holds.
\begin{lemma}\label{lem:subgradient_bound_TCC}
Let $\xbs,\widehat\xbs\in\CB_\nu(\xref)$ with $F(\widehat\xbs)=\ybs$, and let condition \eqref{eqn:TCC} hold. Then
\begin{equation*}
    \DP{\partial\Psi(\xbs)}{\xbs-\widehat\xbs} \ge q(1-\gamma)\Psi(\xbs).
\end{equation*}
\end{lemma}
\begin{proof}
For any fixed $i\in\{1,\dots,N\}$, the Cauchy-Schwarz inequality, the definition of duality map and the tangential cone condition \eqref{eqn:TCC} yield
\begin{align}
\DP{\partial \Psi_i(\xbs)}{\xbs-\widehat\xbs}
    =&\DP{\svaldmap{q}(F_i(\xbs)-\ybs_i)}{F_i(\xbs)-\ybs_i}-\DP{\svaldmap{q}(F_i(\xbs)-\ybs_i)}{F_i(\xbs)-\ybs_i-F_i'(\xbs)(\xbs-\widehat\xbs)}\nonumber\\
    \ge&\yNi{F_i(\xbs)-\ybs_i}^q-\ysNi{\svaldmap{q}(F_i(\xbs)-\ybs_i)}\yNi{F_i(\xbs)-\ybs_i-F_i'(\xbs)(\xbs-\widehat\xbs)}\nonumber\\
    \geq &\yNi{F_i(\xbs)-\ybs_i}^q-\gamma\yNi{F_i(\xbs)-\ybs_i}^{q-1}\yNi{F_i(\xbs)-\ybs_i}=(1-\gamma)\yNi{F_i(\xbs)-\ybs_i}^q.\label{eqn:psi_quasiconvex}
\end{align}
It then follows from \eqref{eqn:obj_function} that
$\DP{\partial\Psi(\xbs)}{\xbs-\widehat\xbs}=
    \frac{1}{N}\sum_{i=1}^N \DP{\partial \Psi_i(\xbs)}{\xbs-\widehat\xbs}\geq q(1-\gamma)\Psi(\xbs)$.
\end{proof}

In order to prove that with suitable step-sizes the iterates $(\xbs_k)_{k\in\bbN}$ stay within a neighborhood around the solution, we first need a descent property with respect to the Bregman distance, which follows from standard arguments, cf. \cite[Lemma 3.6]{JinKereta:2023}.

\begin{lemma}\label{lem:descent_property}
Let Assumption \ref{assn:X_and_Y} hold. Then for any $\widehat\xbs\in\CX$, the iterates in \eqref{eqn:sgd} satisfy
\begin{align}\label{eqn:descent_property}
\bregman{\xbs_{k+1}}{\widehat\xbs}
\leq\bregman{\xbs_k}{\widehat\xbs}-\mu_{k+1}\DP{\gbs_{k+1}}{\xbs_k-\widehat\xbs} + \frac{G_{p^\ast}}{p^\ast}\mu_{k+1}^{p^\ast}\xsN{\gbs_{k+1}}^{p^\ast}.
\end{align}
\end{lemma}

The next result allows using smoothness properties of the forward map $F$ on all iterates. Throughout, let $\Delta_k:=\bregman{\xbs_k}{\xref}$.
\begin{lemma}\label{lem:bregman_monotonicity_qeqp}
Let Assumptions \ref{assn:X_and_Y} and \ref{assn:F} hold in $\CB_\nu(\xref)$, with $\xbs_0\in\CB_\nu(\xref)$, and let step-sizes $(\mu_k)_{k\in\mathbb{N}}$ satisfy $1-\gamma - \Lmax^{p^\ast}\frac{G_{p^\ast}}{p^\ast}\mu_k^{p^\ast-1}\geq C>0$ for all $k\in\bbN$. Then $\xbs_k\in\CB_\nu(\xref)$ for all $k\geq1$.
\end{lemma}
\begin{proof}
Take $k\geq0$ and assume $\xbs_k\in\CB_\nu(\xref)$.
The descent property \eqref{eqn:descent_property} (with $\widehat \xbs = \xref$) yields
\begin{align}\label{eqn:bregman_descent}
\Delta_{k+1}
&\leq\Delta_k-\mu_{k+1}\DP{\gbs_{k+1}}{\xbs_k-\xref} + \frac{G_{p^\ast}}{p^\ast}\mu_{k+1}^{p^\ast}\xsN{\gbs_{k+1}}^{p^\ast}\notag\\
&\leq\Delta_k -p(1-\gamma)\mu_{k+1} \Psi_{i_{k+1}}(\xbs_k) + p\Lmax^{p^\ast}\frac{G_{p^\ast}}{p^\ast}\mu_{k+1}^{p^\ast}\Psi_{i_{k+1}}(\xbs_k)\notag\\
&\leq\Delta_k -p\Big(1-\gamma -\Lmax^{p^\ast}\frac{G_{p^\ast}}{p^\ast}\mu_{k+1}^{p^\ast-1}\Big)\mu_{k+1}\Psi_{i_{k+1}}(\xbs_k),
\end{align}
where we use \eqref{eqn:psi_quasiconvex}, the bound \eqref{eqn:q_var_basic_bound} (which holds since $\xbs_k\in\CB_\nu(\xref)$), and the identity $p^\ast(p-1)=p$.
Since $1-\gamma -\Lmax^{p^\ast}\frac{G_{p^\ast}}{p^\ast}\mu_{k+1}^{p^\ast-1}>0$ by assumption, we have $\Delta_{k+1}\leq\Delta_k$. Thus, by induction, it follows that $\Delta_k\leq\Delta_0$, i.e., $\xbs_k\in\CB_\nu(\xref)$ for all $k\in\bbN$.
\end{proof}

We are now ready to show convergence of the iterations. The lengthy proof is given in Appendix \ref{app:cone_convergence_qeqp}.
\begin{theorem}\label{thm:cone_convergence_qeqp}
Let the conditions of Lemma \ref{lem:bregman_monotonicity_qeqp} hold.
If $\sum_{k=1}^\infty\mu_k=\infty$, then the sequence $(\xbs_k)_{k\in\bbN}$ converges a.s. to a minimizer of $\Psi(\xbs)$, that is,
\begin{equation}\label{eqn:as_convergence_exact}
\bbP\Big(\lim_{k\rightarrow\infty} \inf_{\widetilde \xbs\in \CX_{\min}}\xN{\xbs_k-\widetilde \xbs}=0\Big)=1.
\end{equation}
Moreover, if $\goodset$ is closed and convex, $F$ is weak-to-weak continuous, and  $\Null(F^\prime(\xref))\subset\Null(F^\prime(\xbs))$ holds for all $\xbs\in\goodset$, then the sequence $(\xbs_k)_{k\in\bbN}$ converges a.s. to the $\xbs_0$-minimum-distance solution $\xref$, i.e., $\lim_{k\to\infty}\bregman{\xbs_k}{\xref}=0$ a.s.
\end{theorem}

\begin{remark}\label{rem:convergence_thm}
The assumptions and conclusions in Theorem \ref{thm:cone_convergence_qeqp} can be decomposed into two groups. First, the conditions on $(\mu_k)_{k\in\mathbb{N}}$ are needed for establishing the a.s. convergence of $(\xbs_k)_{k\in\mathbb N}$ in either Bregman distance or in the norm to some nondeterministic limit $\widehat{\xbs}$. The remaining assumptions, on set $\goodset$ and the operator $F$, are necessary to identify this limit as the $\xbs_0$-minimum-distance solution $\xref$; see \cite[Theorem 8]{M18} for a similar analysis for the inexact Newton method.

We also note that Theorem \ref{thm:cone_convergence_qeqp} allows for constant step-sizes, e.g., $\mu_k=\bar{C}(\frac{1-\gamma}{2}\frac{p^\ast}{\Lmax^{p^\ast}G_{p^\ast}})$ for all $k\in\bbN$, with $\bar{C}<1$, which directly yields $\lim_{k\rightarrow\infty}\bbE[\Psi(\xbs_k)]=0$ by taking the full expectation on \eqref{eqn:bregman_descent}. Then, there exists a subsequence $(\xbs_{n_k})_{k\in\bbN}$ such that $\lim_{k\rightarrow\infty} \Psi(\xbs_{n_k})=0$ a.s. Similarly, we then obtain the convergence of $(\xbs_k)_{k\in\bbN}$ to the $\xbs_0$-minimum-distance solution $\xref$ in expectation
\[\lim_{k\to\infty}\bbE[\bregman{\xbs_{k}}{\xref}]=0\quad \mbox{and}\quad \lim_{k\to\infty}\bbE[\xN{\xbs_k-\xref}^p]=0.\]
\end{remark}

Next we extend Theorem \ref{thm:cone_convergence_qeqp} to convergence in expectation.  The proof is in Appendix \ref{app:L1_cone_convergence_qeqp}.
\begin{theorem}\label{thm:L1_cone_convergence_qeqp}
Let the conditions of Theorem \ref{thm:cone_convergence_qeqp} hold with a closed and convex $\goodset$, a weak-to-weak continuous $F$ and $\Null(F^\prime(\xref))\subset\Null(F^\prime(\xbs))$ for all $\xbs\in\goodset$.
Then there holds
$$\lim_{k\rightarrow\infty}\bbE[\bregman{\xbs_k}{\xref}]=0.$$
Moreover, for $1\leq r\leq p$, we have $$\lim_{k\rightarrow\infty} \bbE[\xN{\xbs_k-\xref}^r]=0\quad \mbox{and}\quad \lim_{k\rightarrow\infty} \bbE[\xsN{\dmapX{p}(\xbs_k)-\dmapX{p}(\xref)}^{p^\ast}]=0.$$
\end{theorem}

\subsection{Regularizing property for noisy data}\label{sec:reg_property}
Now we show that SGD has a regularizing effect, in the sense that the expected error $\bbE[\bregman{\xbs_{k(\delta)}^\delta}{\xref}]$ converges to $0$ as the noise level $\delta$ goes to $0$, with properly selected stopping indices $k(\delta)$.
Let $(\xbs_k)_{k\in\bbN}$ and $(\xbs_k^\delta)_{k\in\bbN}$ be the SGD iterates for noise-free and noisy measurements, defined, respectively, by
\begin{align}
&\xbs_{k+1} = \dmapXs{p}{\left(\dmapX{p}(\xbs_k) - \mu_{k+1}\gbs_{k+1}\right)}, \quad {\rm with }\quad \gbs_{k+1} = \gbs(\xbs_k,\ybs,i_{k+1}),\label{eqn:sgd_cleaniterates}\\
&\xbs_{k+1}^\delta = \dmapXs{p}{\left(\dmapX{p}(\xbs_k^\delta) - \mu_{k+1}\gbs_{k+1}^\delta\right)}, \quad {\rm with } \quad \gbs_{k+1}^\delta = \gbs(\xbs_k^\delta,\ybs^\delta,i_{k+1})\label{eqn:sgd_noisyiterates}.
\end{align}

We first we show that all iterates $(\xbs_k^\delta)_{0<k\leq k(\delta)}$ generated by \eqref{eqn:sgd_noisyiterates} stay within the closed ball $\CB_\nu(\xref)$, for a suitably defined radius $\nu>0$.
This does not need to be shown for linear forward operators, as they are globally well-behaved.
However, for nonlinear operators, this condition is crucial, since the regularization analysis requires all iterates to stay within such a ball in order to satisfy Assumption \ref{assn:F} (i.e., so that the tangential cone and Lipschitz boundedness can be invoked).
This is a key difference from the linear case.
Below we denote $\Delta_k^\delta:=\bregman{\xbs_k^\delta}{\xref}$.

\begin{lemma}\label{lem:bregman_monotonicity_qeqp_noisy}
Let Assumption \ref{assn:X_and_Y} hold, and $\xbs_0\in\CB_\nu(\xref)$.
Let Assumption \ref{assn:F} hold with $0<\gamma<1/2$ and $\nu = \bregman{\xbs_0}{\xref} + \frac{\omega^{-p}}{p}(1+\gamma)^p\Gamma$, for some small $\omega>0$.
Assume that the step-sizes $(\mu_\ell)_{\ell=1}^{ k(\delta)}$ satisfy $$\delta^p\sum_{\ell=1}^{k(\delta)}\mu_\ell\le \Gamma\quad \mbox{and}\quad 1-\gamma -\Lmax^{p^\ast}\frac{G_{p^\ast}}{p^\ast}\mu_{\ell}^{p^\ast-1}-\frac{p-1}{p^2}\omega^{p^\ast} \ge \bar{C}_{ub} > 0.$$ Then $\xbs_k^\delta\in\CB_\nu(\xref)$ holds for all $0<k\leq k(\delta)$.
\end{lemma}
\begin{proof}
The proof proceeds by mathematical induction. The assertion holds trivially in the base case $m=0$.
Assume that the assertion $\xbs_m^\delta\in\CB_\nu(\xref)$ holds for all $0<m\le k$. For noisy iterates $(\xbs_k^\delta)_{k=1}^{k(\delta)}$, the descent property \eqref{eqn:descent_property} with $\widehat\xbs = \xref$ still holds. We then obtain
\begin{align}\label{eqn:bregman_descent_noisy}
\Delta_{k+1}^\delta
\leq &\Delta_k^\delta-\mu_{k+1}\DP{\gbs_{k+1}^\delta}{\xbs_k^\delta-\xref} + \frac{G_{p^\ast}}{p^\ast}\mu_{k+1}^{p^\ast}\xsN{\gbs_{k+1}^\delta}^{p^\ast}\notag\\
\leq& \Delta_k^\delta-p(1-\gamma)\mu_{k+1}\Psi_{i_{k+1}}(\xbs_k^\delta) + p\Lmax^{p^*}\frac{G_{p^\ast}}{p^\ast}\mu_{k+1}^{p^\ast}\Psi_{i_{k+1}}(\xbs_k^\delta) \notag\\
&+ (1+\gamma)\delta\mu_{k+1}\yNi{F_{i_{k+1}}(\xbs_{k}^\delta)-\ybs_{i_{k+1}}^\delta}^{p-1}\notag\\
\leq& \Delta_k^\delta - p\Big(1-\gamma -\Lmax^{p^\ast}\frac{G_{p^\ast}}{p^\ast}\mu_{k+1}^{p^\ast-1}\Big)\mu_{k+1}\Psi_{i_{k+1}}(\xbs_k^\delta) \nonumber\\
&+ (1+\gamma)\delta\mu_{k+1}\yNi{F_{i_{k+1}}(\xbs_{k}^\delta)-\ybs_{i_{k+1}}^\delta}^{p-1},
\end{align}
where we have used the estimate \eqref{eqn:psi_quasiconvex}, the bound \eqref{eqn:q_var_basic_bound} (which holds since $\xbs_k^\delta\in\CB_\nu(\xref)$) and $\yNi{\ybs_i-\ybs_i^\delta}\le\delta_i\le \delta$, and the identity $p^\ast(p-1)=p$.
Using Young's inequality $ab\leq \frac{\omega^{-p}}{p}a^p+\frac{\omega^{P^\ast}}{p^\ast}b^{p^\ast}$ with $a=(1+\gamma)\delta$ and $b=\yNi{F_{i_{k+1}}(\xbs_{k}^\delta)-\ybs_{i_{k+1}}^\delta}^{p-1}$ and $\omega>0$, we have
\begin{align}\label{eqn:young_ineq_residual}
    (1+\gamma)\delta\yNi{F_{i_{k+1}}(\xbs_{k}^\delta)-\ybs_{i_{k+1}}^\delta}^{p-1}
    \le \frac{\omega^{-p}}{p}(1+\gamma)^p\delta^p+\frac{p-1}{p}\omega^{p^\ast}\yNi{F_{i_{k+1}}(\xbs_{k}^\delta)-\ybs_{i_{k+1}}^\delta}^{p}.
\end{align}
By plugging \eqref{eqn:young_ineq_residual}
into \eqref{eqn:bregman_descent_noisy}, and using the assumptions $1-\gamma -\Lmax^{p^\ast}\frac{G_{p^\ast}}{p^\ast}\mu_{\ell}^{p^\ast-1}-\frac{p-1}{p^2}\omega^{p^\ast}\geq \bar{C}_{ub}>0 $ and $\delta^p\sum_{\ell=1}^{k(\delta)}\mu_\ell\le \Gamma$, we deduce
\begin{align*}
\Delta_{k+1}^\delta
&\le \Delta_k^\delta - p\Big(1-\gamma -\Lmax^{p^\ast}\frac{G_{p^\ast}}{p^\ast}\mu_{k+1}^{p^\ast-1}-\frac{p-1}{p^2}\omega^{p^\ast}\Big)\mu_{k+1}\Psi_{i_{k+1}}(\xbs_k^\delta) + \frac{\omega^{-p}}{p}(1+\gamma)^p\delta^p\mu_{k+1}\\
&\le \Delta_k^\delta + \frac{\omega^{-p}}{p}(1+\gamma)^p\delta^p\mu_{k+1}
\le \Delta_0 + \frac{\omega^{-p}}{p}(1+\gamma)^p\delta^p\sum_{\ell=1}^{k(\delta)}\mu_\ell
\le \Delta_0 + \frac{\omega^{-p}}{p}(1+\gamma)^p\Gamma.
\end{align*}
It follows that $\Delta_{k+1}^\delta\le\nu$. Thus, $\xbs_{k}^\delta\in\CB_\nu(\xref)$ for all $0<k\le k(\delta)$ with $\nu = \Delta_0 + \frac{\omega^{-p}}{p}(1+\gamma)^p\Gamma$.
\end{proof}

\begin{remark}
In Lemma \ref{lem:bregman_monotonicity_qeqp_noisy}, the step-sizes $(\mu_\ell)_{\ell=1}^{ k(\delta)}$ must satisfy
$$1-\gamma -\Lmax^{p^\ast}\frac{G_{p^\ast}}{p^\ast}\mu_{\ell}^{p^\ast-1}-\frac{p-1}{p^2}\omega^{p^\ast}>0,$$ i.e., $\mu_{\ell}^{p^\ast-1}<\frac{p^\ast}{G_{p^\ast}\Lmax^{p^\ast}}(1-\gamma -\frac{p-1}{p^2}\omega^{p^\ast})$.
Since the step-sizes $(\mu_\ell)_{\ell\in\mathbb{N}}$ must be positive, we need to ensure that $1-\gamma -\frac{p-1}{p^2}\omega^{p^\ast}>0$ holds.
By the tangential cone condition in Assumption \ref{assn:F} {\rm(iii)}, we also have $\gamma>0$, which implies $\omega^{p^\ast}<\frac{p^2}{p-1}$.
By writing now $\omega^{p^\ast} = \epsilon\frac{p^2}{p-1}$, for some $\epsilon < 1$, we find that $\gamma < 1 - \frac{p-1}{p^2}\omega^{p^\ast} = 1 - \epsilon < 1$. However, note that the radius $\nu= \bregman{\xbs_0}{\xref} + \frac{\omega^{-p}}{p}(1+\gamma)^p\Gamma$ increases with the decrease of $\omega$ {\rm(}i.e., as $\epsilon$ decreases{\rm)}.
In the limiting case as $\epsilon\to 0$, this would imply that Assumption \ref{assn:F} needs to hold globally on $\CX$.
In the statement of Lemma \ref{lem:bregman_monotonicity_qeqp_noisy},  we take a convenient choice of $\epsilon=\frac{1}{2}$, which leads to the admissible range $\gamma \in (0,1/2)$. This condition is stronger than that in Assumption \ref{assn:F} {\rm(iii)}, though it can be made arbitrarily close to it by adjusting the constant or, if necessary, reducing the step-sizes to decrease $\Gamma$. However, such reductions may lead to slower convergence.
Nonetheless, restricting $\gamma$ to a specific admissible range is a standard approach in the literature \cite[Section 2]{KaltenbacherNeubauerScherzer:2008}.
\end{remark}

Next we show the stability of SGD with respect to noise. This requires that $\CY$ is uniformly smooth, for which the corresponding duality maps are single valued and continuous. The proof follows the argument of \cite[Lemma 4.2]{JinKereta:2023}, using also Assumption \ref{assn:F}(iii)
\begin{lemma}\label{lem:coupled_noise_convergence}
Let Assumptions \ref{assn:X_and_Y} and \ref{assn:F} hold. Consider the iterations \eqref{eqn:sgd_cleaniterates} and \eqref{eqn:sgd_noisyiterates} with the same initialization $\xbs_0^\delta=\xbs_0$, and following the same path {\rm(}i.e., using same random indices $i_{k}${\rm)}. Then, for any fixed $k\in\bbN$, we have
\begin{align*}
    \lim_{\delta\searrow0}\bbE[\bregman{\xbs^{\delta}_k}{\xbs_k}]
    =\lim_{\delta\searrow0}\bbE[\xN{\xbs^{\delta}_k -\xbs_k}]
    =\lim_{\delta\searrow0}\bbE[\xsN{\dmapX{p}(\xbs_k^\delta)-\dmapX{p}(\xbs_k)}]=0.
\end{align*}
\end{lemma}

Now we show the regularizing property of SGD for suitable stopping indices $k(\delta)$.
\begin{theorem}\label{thm:regularisation_property}
Let the conditions in Lemma \ref{lem:bregman_monotonicity_qeqp_noisy} hold, and $(\mu_k)_{k\in\bbN}$ satisfy $\sum_{k=1}^\infty\mu_{k}=\infty$. If
$
\lim_{\delta\searrow0}k(\delta)=\infty$ and $
\lim_{\delta\searrow0} \delta^p\sum_{\ell=1}^{k(\delta)}\mu_\ell=0$,
and if the set $\goodset$ is closed and convex, $F$ is weak-to-weak continuous, and $\Null(F^\prime(\xref))\subset\Null(F^\prime(\xbs))$ for all $\xbs\in\goodset$, then
\begin{align*}
    \lim_{\delta\searrow0} \bbE[\bregman{\xbs_{k(\delta)}^{\delta}}{\xref}]=0.
\end{align*}
\end{theorem}
\begin{proof}
Take any $\delta>0$ and $k\in\bbN$.
Using the three-point identity \eqref{eqn:3_point_id} of the Bregman distance, and the Cauchy-Schwarz inequality, we have
\begin{align}\label{eqn:noisy_3point}
    \Delta_k^\delta
    &= \bregman{\xbs_k^\delta}{\xbs_k} + \Delta_k  + \DP{\dmapX{p}(\xbs_k)-\dmapX{p}(\xbs_k^\delta)}{\xbs_k-\xref}\nonumber \\
    &\leq \bregman{\xbs_k^\delta}{\xbs_k} + \Delta_k + \xsN{\dmapX{p}(\xbs_k)-\dmapX{p}(\xbs_k^\delta)}\xN{\xbs_k-\xref}.
\end{align}
Consider a sequence $(\delta_j)_{j\in\bbN}$ decaying to zero.
For any $\epsilon>0$, it suffices to find $j_\epsilon\in\bbN$ such that $\bbE[\Delta_{k(\delta_j)}^{\delta_j}]\leq 4\epsilon$ for all $j\geq j_\epsilon$. By Theorem \ref{thm:L1_cone_convergence_qeqp},
there exists $k_\epsilon\in\bbN$ such that for all $k\geq k_\epsilon$, we have
\begin{align}\label{eqn:noiseless_kbound}
    \bbE[\Delta_k]<\epsilon\quad \text{and}\quad \bbE[\xN{\xbs_k-\xref}]<\epsilon^{1/2}.
\end{align}
Moreover, for any fixed $k_\epsilon$, by Lemma \ref{lem:coupled_noise_convergence}, there exists $j_1\in\bbN$ such that for all $j\geq j_1$, we have
\begin{align}\label{eqn:noisy_kepsilon_bound}
    \bbE[\bregman{\xbs_{k_\epsilon}^{\delta_j}}{\xbs_{k_\epsilon}}] <\epsilon
    \quad \text{and}\quad
    \bbE[\xsN{\dmapX{p}(\xbs_{k_\epsilon})-\dmapX{p}(\xbs_{k_\epsilon}^{\delta_j})}]<\epsilon^{1/2}.
\end{align}
Thus, by plugging the estimates \eqref{eqn:noiseless_kbound} and \eqref{eqn:noisy_kepsilon_bound} into \eqref{eqn:noisy_3point}, we have $\bbE[\Delta_{k_\epsilon}^{\delta_j}] < 3\epsilon$, for all $j\geq j_1$.
Note, however, that the same does not necessarily hold for all $k\geq k_\epsilon$, and thus for a monotonically increasing sequence of stopping indices $k(\delta_j)$, the expectation sequences $\bbE[\Delta_{k(\delta_j)}^{\delta_j}]$ are not necessarily monotone.
Instead, from Lemma \ref{lem:bregman_monotonicity_qeqp_noisy}, the descent property with noisy data yields
\begin{equation*}
    \Delta_{k+1}^\delta
    \le \Delta_k^\delta - p\bar{C}_{ub}\mu_{k+1}\Psi_{i_{k+1}}(\xbs_k^\delta) + \frac{\omega^{-p}}{p}(1+\gamma)^p\delta^p\mu_{k+1},
\end{equation*}
with $\bar{C}_{ub}>0$ (cf. condition of Lemma \ref{lem:bregman_monotonicity_qeqp_noisy}).
Taking the full expectation of the inequality gives
\begin{align}\label{eqn:bregman_descent_noisy_expectation}
\bbE[\Delta_{k+1}^\delta]
\leq \bbE[\Delta_k^\delta]-p\bar{C}_{ub}\mu_{k+1}\bbE[\Psi(\xbs_k^\delta)] + \frac{\omega^{-p}}{p}(1+\gamma)^p\delta^p\mu_{k+1}.
\end{align}
Replacing $k+1$ with $k(\delta)$, and using the inductive argument, implies
\begin{align*}
    \bbE[\Delta^\delta_{k(\delta)}]
    \leq \bbE[\Delta^\delta_{k(\delta)-1}] + \frac{\omega^{-p}}{p}(1+\gamma)^p\delta^p\mu_{k(\delta)}
    \leq\bbE[\Delta^\delta_{k_\epsilon}] + \frac{\omega^{-p}}{p}(1+\gamma)^p\delta^p\sum_{\ell=1}^{k(\delta)}\mu_{\ell}.
\end{align*}
Under the conditions on $k(\delta)$, there exists $j_2\in\bbN$ such that for all $j\geq j_2$, we have $k(\delta_j)\geq k_\epsilon$ and $\frac{\omega^{-p}}{p}(1+\gamma)^p\delta_j^p\sum_{\ell=1}^{k(\delta_j)}\mu_{\ell} <\epsilon$.
Taking ${j_\epsilon}=j_1\vee j_2$ gives $\bbE[\Delta_{k(\delta_j)}^{\delta_j}] < 4\epsilon$ for all $j\geq j_\epsilon$.
Hence, the desired claim follows.
\end{proof}

\begin{remark}
By Theorem \ref{thm:regularisation_property}, the condition $
\lim_{\delta\searrow0} \delta^p\sum_{\ell=1}^{k(\delta)}\mu_\ell=0$ is required. For a constant step-size $\mu_\ell = \mu_0$, this implies $k(\delta) = o(\delta^{-p})$.
\end{remark}

\section{Convergence rates under conditional stability}\label{sec:nonlin_rate}

Theorems \ref{thm:L1_cone_convergence_qeqp} and \ref{thm:regularisation_property} establish convergence in expectation for exact and noisy data, respectively, but do not provide convergence rates. To derive convergence rates, additional assumptions on the solution $\xref$ are required, which are collectively known as source conditions. One approach is via conditional stability that is known for many inverse problems for PDEs (see, e.g., \cite{BellassouedYamamoto:2017,KlibanovTimonov:2004,Yamamoto:2009}), and which has been used to analyze both variational regularization \cite{ChengYamamoto:2000,WernerHofmann:2020} and iterative regularization \cite{deHoopQiuScherzer:2012}.
Typically, a H\"{o}lder-type conditional stability estimate  \cite{deHoopQiuScherzer:2012} assumes the existence of constants $\alpha\geq1$ and $C_\alpha>0$ such that
\begin{align}\label{eqn:Fcond_stability}
\bregman{\xbs}{\tilde\xbs}^\alpha\le C_\alpha^{-1} \yN{F(\xbs)-F(\tilde\xbs)}^p, \quad\forall \xbs,\tilde\xbs\in\goodset.
\end{align}
Note that under \eqref{eqn:Fcond_stability}, the $\xbs_0$-minimum-distance solution is equivalent to the $\xbs_0$-minimum-norm solution and is unique.
We denote it by $\xref$.

\begin{theorem}\label{thm:cond_stability_convergence_qeqp}
Let \eqref{eqn:Fcond_stability}, and Assumptions \ref{assn:X_and_Y} and \ref{assn:F} hold.
Let the step-sizes $(\mu_k)_{k\in\bbN}$ be such that $\frac{C_pC_\alpha}{N}\Big(1-\gamma - \Lmax^{p^\ast}\frac{G_{p^\ast}}{p^\ast}\mu_k^{p^\ast-1}\Big)\ge C_{ub}>0$, with $C_p\ge1$ given by Remark \ref{rmk:product}, and $\sum_{k=1}^\infty \mu_k=\infty$.
Then $\lim_{k\rightarrow\infty}\bbE[\bregman{\xbs_{k}}{\xref}]=0$ holds if $\xbs_0\in\CB_\nu(\xref)$.
Moreover, if $\alpha=1$ and $\mu_j\geq c_{lb}/C_{ub}>0$ for $j\in\bbN$ then $\bbE[\bregman{\xbs_{k}}{\xref}]\leq (1-c_{lb})^k\bregman{\xbs_{0}}{\xref}$, and otherwise
\begin{align*}
\bbE[\bregman{\xbs_{k}}{\xref}]\leq \left\{
\begin{aligned}\Big(1+C_{ub}(\alpha-1)\bregman{\xbs_{0}}{\xref}^{\alpha-1}\sum_{j=1}^k \mu_j\Big)^{\frac{1}{1-\alpha}}\bregman{\xbs_{0}}{\xref},&\quad \text{ if } \alpha>1,\\
\exp\Big(-\sum_{j=1}^{k} C_{ub}\mu_j\Big)\bregman{\xbs_{0}}{\xref},&\quad \text{ if } \alpha=1.
\end{aligned}\right.
\end{align*}
\end{theorem}
\begin{proof}
Taking the full expectation on the inequality \eqref{eqn:cond_expect_bregman}, we get
\begin{align*}
\bbE[\Delta_{k+1}]\leq\bbE[\Delta_k] -p \Big(1-\gamma - \Lmax^{p^\ast}\frac{G_{p^\ast}}{p^\ast}\mu_k^{p^\ast-1}\Big)\mu_{k+1}\bbE[\Psi(\xbs_k)].
\end{align*}
The inequality $\Psi(\xbs)\geq \frac{C_p}{pN}\yN{F(\xbs)-\ybs}^p$,  conditional stability \eqref{eqn:Fcond_stability} and Jensen's inequality yield
\begin{align*}
\bbE[\Delta_{k+1}]
&\leq\bbE[\Delta_k]-p\Big(1-\gamma - \Lmax^{p^\ast}\frac{G_{p^\ast}}{p^\ast}\mu_k^{p^\ast 1}\Big)\mu_{k+1}\bbE[\Psi(\xbs_k)]\\
&\leq\bbE[\Delta_k]-C_{ub}\mu_{k+1}\bbE[\Delta_k^\alpha]
 \leq\bbE[\Delta_k]-C_{ub}\mu_{k+1}\bbE[\Delta_k]^\alpha.
\end{align*}
Since $C_{ub}>0$, the sequence $(\bbE[\Delta_k])_{k\in\bbN}$ is monotonically decreasing, and since it is bounded from below (by $0$), it is convergent.
Assume now that $\lim_{k\rightarrow\infty}\bbE[\Delta_k]=L>0$,
Then $\bbE[\Delta_k]\geq L\geq0$, since the sequence is monotonic, which for any $k_0\in\bbN$ gives
\begin{align*}
\bbE[\Delta_{k+1}]
&\leq \bbE[\Delta_k]-C_{ub}\mu_{k+1}\bbE[\Delta_k]^{\alpha}\le \bbE[\Delta_k]-C_{ub}\mu_{k+1}L^\alpha\\
&\le\bbE[\Delta_{k_0}]-C_{ub}L^\alpha\sum_{j=k_0}^{k+1} \mu_{j}\rightarrow -\infty,
\end{align*}
since $\sum_{j=1}^\infty \mu_{j}=\infty$ by assumption, which gives a contradiction.
Therefore, $\lim_{k\rightarrow\infty}\bbE[\Delta_k]=0$.
For $\alpha=1$, direct computation yields
\[ \bbE[\Delta_{k+1}] \leq (1-C_{ub}\mu_{k+1})\bbE[\Delta_{k}]\leq \prod_{j=1}^{k+1} (1-C_{ub}\mu_j) \Delta_0.\]
If $(\mu_k)_{k\in\mathbb{N}}$ are constant (or bounded from below), we thus get $\bbE[\Delta_{k}]\leq (1-c_{lb})^k\Delta_{0}$.
Otherwise, by the elementary inequality $1-x\leq e^{-x}$ for $x\geq0$, we have
\[ \bbE[\Delta_{k+1}] \leq\prod_{j=1}^{k+1} (1-C_{ub}\mu_j) \Delta_0 \leq \exp\Big(-\sum_{j=1}^{k+1} C_{ub}\mu_j\Big)\Delta_0.\]
For $\alpha>1$, by Polyak's inequality given in Lemma \ref{lem:polyak_series}, we have
\[ \bbE[\Delta_{k+1}]\leq \Delta_0\Big(1+C_{ub}(\alpha-1)\Delta_0^{\alpha-1} \sum_{j=1}^{k+1} \mu_{j}\Big)^{\frac{1}{1-\alpha}}.\]
This completes the proof of the theorem.
\end{proof}

\begin{remark}
Theorem \ref{thm:cond_stability_convergence_qeqp} allows for constant step-sizes.
Choosing $\mu_k=\mu_0$ such that $\mu_0^{p^\ast-1} = \frac{p^\ast(1-\gamma)}{2\Lmax^{p^\ast}G_{p^\ast}}$ ensures $1-\gamma -\Lmax^{p^\ast}\frac{G_{p^\ast}}{p^\ast}\mu_{k+1}^{p^\ast-1}= \frac{1-\gamma}{2}>0$.
When $\alpha=1$, we get a linear convergence rate.
Namely, with $C_{ub}=\frac{C_pC_\alpha(1-\gamma)}{2N}$ and $c_{lb}=\mu_0/C_{ub}$, we get
\begin{align*}
    &\bbE[\bregman{\xbs_{k+1}}{\xref}]
    \leq \Big(1- \frac{C_pC_\alpha(1-\gamma)}{2N}\mu_0\Big)\bbE[\bregman{\xbs_{k}}{\xref}]\\
    \leq &\bigg(1- \frac{C_pC_\alpha}{N}\bigg(\frac{1-\gamma}{2}\bigg)^{1+1/(p^\ast-1)}\Lmax^{-p^\ast/(p^\ast-1)}\bigg(\frac{p^\ast}{G_{p^\ast}}\bigg)^{1/(p^\ast-1)}\bigg)\bbE[\bregman{\xbs_{k}}{\xref}]\\
    \leq&\bigg(1- \frac{C_pC_\alpha}{N}\bigg(\frac{1-\gamma}{2}\bigg)^p\Lmax^{-p}\bigg(\frac{p^\ast}{G_{p^\ast}}\bigg)^{p/p^\ast}\bigg)^{k+1}\bregman{\xbs_{0}}{\xref}.
\end{align*}
For $\alpha>1$ one can instead observe a phase transition.
To observe this, let $k_L$ be the largest positive integer such that $\bbE[\Delta_{k_L-1}]>1$. A direct computation then yields
$\bbE[\Delta_{k}]\leq (1-c_{lb})^{(k-1)\alpha+1}\Delta_0$, for all $k\leq k_L$. Thereafter the convergence rate deteriorates, and we can only have a sub-linear rate:
\begin{align*} \bbE[\Delta_{k+1}]\leq\bbE[\Delta_{k}]\Big(1-c_{lb}\bbE[\Delta_{k}]^{\alpha-1}\Big) \Rightarrow \bbE[\Delta_{k+1}]^{1-\alpha}\geq\bbE[\Delta_{k}]^{1-\alpha}\Big(1-c_{lb}\bbE[\Delta_{k}]^{\alpha-1}\Big)^{1-\alpha}.
\end{align*}
Assume now that $c_{lb}\geq1$.
We then have $(1-c_{lb}\bbE[\Delta_{k}]^{\alpha-1})^{1-\alpha}\geq 1-c_{lb}(1-\alpha)\bbE[\Delta_{k}]^{\alpha-1}$, giving
\begin{align*}
    \bbE[\Delta_{k+1}]^{1-\alpha}
    \ge\bbE[\Delta_{k}]^{1-\alpha}+c_{lb}(\alpha-1)
    \ge\bbE[\Delta_{k_L}]^{1-\alpha}+(k-k_L)c_{lb}(\alpha-1),
\end{align*}
and thus $\bbE[\Delta_{k+1}]\leq\big(\bbE[\Delta_{k_L}]^{1-\alpha}+(k-k_L)c_{lb}(\alpha-1)\big)^{1/1-\alpha}$, which agrees with the rate obtained through Polyak's inequality (cf. Lemma \ref{lem:polyak_series}).
\end{remark}

The next result gives the convergence rate for noisy data under  conditional stability.
\begin{theorem}\label{thm:noisy_convergence_rate}
Let the conditions in Lemma \ref{lem:bregman_monotonicity_qeqp_noisy} hold, and let $F$ satisfy \eqref{eqn:Fcond_stability}.
If the step-sizes $(\mu_k)_{k\in\bbN}$ are chosen such that
$\sum_{k=1}^\infty \mu_k=\infty$, then for $\alpha\ge1$, there holds
\begin{equation}\label{eqn:cvg_rate_noisy}
\bbE[\bregman{\xbs_{k(\delta)}^\delta}{\xref}]
\leq 2\bar{C}_{p\gamma}^\alpha{(\bar{C}_{ub}^{\ast})^{-\alpha}}\delta^{p/\alpha},
\end{equation}
with $\bar{C}_{p\gamma} = \bar{C}_{ub}+\frac{\omega^{-p}}{p}(1+\gamma)^p$ and $\bar{C}_{ub}^{\ast}=2^{1-p}C_\alpha\bar{C}_{ub}\frac{C_p}{N}$, where $\bar{C}_{ub}>0$ and $C_p\ge1$ are given by Lemma \ref{lem:bregman_monotonicity_qeqp_noisy} and Remark \ref{rmk:product}, respectively.
\end{theorem}
\begin{proof}
Using the estimate $\Psi(\xbs_k^\delta)\ge \frac{C_p}{pN}\yN{F(\xbs_{k}^\delta)-\ybs^\delta}^p$ and \eqref{eqn:bregman_descent_noisy_expectation} yields
\begin{align}\label{eqn:Bregman_descent_noisy_expectation_rate}
\bbE[\Delta_{k+1}^\delta]
&\leq \bbE[\Delta_k^\delta]-\bar{C}_{ub}\frac{C_p}{N}\mu_{k+1}\bbE[\yN{F(\xbs_{k}^\delta)-\ybs^\delta}^p] + \frac{\omega^{-p}}{p}(1+\gamma)^p\delta^p\mu_{k+1}.
\end{align}
The conditional stability \eqref{eqn:Fcond_stability} and the bound $\yN{\ybs^\delta-\ybs}^p\le \frac{N}{C_p}\delta^p$, which follows from the definition of $\delta$ and the norm equivalence \eqref{eqn:lr_norm} for $1<r<\infty$, imply
\begin{align*}
    (\Delta_k^\delta)^\alpha
    \le 2^{p-1} C_\alpha^{-1}
    \Big( \yN{F(\xbs_{k}^\delta)-\ybs^\delta}^p + \frac{N}{C_p}\delta^p \Big).
\end{align*}
Taking expectation on the above inequality and using Jensen's inequality lead to
\begin{align}\label{eqn:Fcond_stability_variant}
    -\bbE[\yN{F(\xbs_{k}^\delta)-\ybs^\delta}^p] \le -2^{1-p}C_\alpha\bbE[\Delta_k^\delta]^\alpha + \frac{N}{C_p}\delta^p.
\end{align}
By plugging \eqref{eqn:Fcond_stability_variant} into \eqref{eqn:Bregman_descent_noisy_expectation_rate}, we get
\begin{align*}
\bbE[\Delta_{k+1}^\delta]
&\leq \bbE[\Delta_k^\delta]-\mu_{k+1}(\bar{C}_{ub}^\ast\bbE[\Delta_k^\delta]^\alpha - \bar{C}_{p\gamma}\delta^p).
\end{align*}
Let
$\delta_\xi
\coloneqq \Big(\frac{\bar{C}_{p\gamma}}{\bar{C}_{ub}^\ast}\Big)^{1/\alpha}\delta^{p/\alpha}$.
Assume now that the statement does not hold, that is, that there exists a $\delta>0$ such that  $\bbE[\Delta_{k}^\delta]> 2\delta_\xi$ holds for all $k$.
A direct computation gives
\begin{equation*}
    \frac{\bar{C}_{p\gamma}\delta^p}{\bar{C}_{ub}^\ast\bbE[\Delta_k^\delta]^\alpha}
    < \frac{\bar{C}_{p\gamma}\delta^p}{\bar{C}_{ub}^\ast(2\delta_\xi)^\alpha} = \frac{1}{2^\alpha},
\end{equation*}
that is,
$\bar{C}_{p\gamma}\delta^p < \tfrac{1}{2^\alpha}\bar{C}_{ub}^\ast\bbE[\Delta_k^\delta]^\alpha.$
This directly leads to
\begin{align}\label{eqn:conv_rates_noisy_bregman_relation}
\bbE[\Delta_{k+1}^\delta]
< \bbE[\Delta_k^\delta] - \big(1-\tfrac{1}{2^\alpha}\big)\bar{C}_{ub}^\ast\mu_{k+1}\bbE[\Delta_k^\delta]^\alpha.
\end{align}
Plugging in $\bbE[\Delta_{k}^\delta]> 2\delta_\xi$ then gives
\[\bbE[\Delta_{k+1}^\delta]
< \bbE[\Delta_k^\delta] - \big(2^\alpha-1\big)\bar{C}_{p\gamma}\mu_{k+1}\delta^p.\]
Hence, by induction we have
\[\bbE[\Delta_{k+1}^\delta]
< \Delta_0^\delta - \big(2^\alpha-1\big)\bar{C}_{p\gamma}\sum_{j=0}^k \mu_{j+1},\]
for all $k$. This leads to a contradiction since $\sum_{j=0}^\infty \mu_{j+1}=\infty$ holds by assumption, which means that for all large enough $k$ the right hand side is negative, but the left hand side is non-negative.
\end{proof}

\begin{remark}
Theorem \ref{thm:noisy_convergence_rate} implies
$
\bbE[\bregman{\xbs_{k(\delta)}^\delta}{\xref}] \leq C\delta^{p/\alpha}$ and  $\bbE[\xN{\xbs_{k(\delta)}^\delta-\xref}] \leq C\delta^\alpha.$
In view of the estimate \eqref{eqn:cvg_rate_noisy}, we have $\lim_{\delta\searrow0}\bbE[\bregman{\xbs_{k(\delta)}^\delta}{\xref}] = 0$, and hence by \eqref{eqn:norm_leq_bregman}, we have $$\lim_{\delta\searrow0}\bbE\big[\xN{\xbs_{k(\delta)}^\delta-\xref}\big]=0.$$
Theorem \ref{thm:noisy_convergence_rate} allows for constant step-sizes. Choosing constant step-sizes $\mu_k=\mu_0$ such that $\mu_0^{p^\ast-1} = \frac{p^*(1-\gamma-(p-1)p^{-2}\omega^{p^\ast})}{2\Lmax^{p^\ast}G_{p^\ast}}$ guarantees
$$\bar C_{ub}=\frac{1-\gamma-(p-1)p^{-2}\omega^{p^\ast}}{2}>0,\quad \bar{C}_{p\gamma} = \bar C_{ub}+\frac{\omega^{-p}}{p}(1+\gamma)^p\quad \mbox{and}\quad \bar{C}_{ub}^{\ast}=2^{1-p}C_\alpha \bar C_{ub}\frac{C_p}{N}.$$
When $\alpha=1$, we have
\begin{align*}
\bbE[\bregman{\xbs_{k(\delta)}^\delta}{\xref}]
&\leq \Big(1-\bar{C}_{ub}^\ast\mu_0\Big)\bbE[\bregman{\xbs_{k(\delta)-1}^\delta}{\xref}] +
      \bar{C}_{p\gamma}\mu_{0}\delta^p\\
&\leq \Big(1-\bar{C}_{ub}^\ast\mu_0\Big)^{k(\delta)}\bregman{\xbs_{0}}{\xref} + \sum_{j=0}^{k(\delta)-1}(1-
      \bar{C}_{ub}^\ast\mu_0)^j\bar{C}_{p\gamma}\mu_0\delta^p.
\end{align*}
\end{remark}

\section{Numerical experiments and discussions}\label{sec:numerics}
In this section, we illustrate the behavior of SGD \eqref{eqn:sgd} with two nonlinear image reconstruction problems, namely, Schlieren tomography and electrical impedance tomography, with either Gaussian or impulsive noise.
The reconstructions for all experiments are obtained using one single stochastic run.

\subsection{Schlieren tomography}\label{sec:schlieren-tomography}
First we investigate the behavior of SGD on Schlieren tomography, which aims to reconstruct the pressure field $\xbs$ from optical measurements.
In a Schlieren optical system, pressure variations within the given medium cause small deflections in the rays of light passing through it.
The resulting light intensity emerging from the tank, and recorded at the detector, is proportional to the square of the line integral of the pressure $\xbs$ along the light path.

Let $\sigma_i = \sigma(\theta_i) = (\cos\theta_i, \sin\theta_i)\in \mathbb{S}^1$, $i=1,\ldots,N$, be a set of recording directions,
where $(\theta_i)_{i=1}^N\subset[0,\pi)$ are the projection angles.
For an open and bounded region $D\subset \bbR^2$, the Schlieren operator in the direction $\sigma_i$ can be expressed as the square of the Radon transform of the pressure field, $F_i(\xbs) \coloneqq R_i^2(\xbs)$.
The Radon transform $R_i: H_0^1(D)\to L^2([-1,1])$ is given by
\begin{equation*}
     R_i(\xbs)(s,\sigma_i) := \int_\bbR \xbs(s\sigma_i+r\sigma_i^\perp){\rm d}r, \quad (s,\sigma_i)\in [-1,1]\times \mathbb{S}^1.
\end{equation*}
This leads to a nonlinear system of equations
\begin{equation}
    F(\xbs):=(F_1(\xbs),\ldots,F_N(\xbs))=\ybs:=(\ybs_1,\ldots,\ybs_N).
\end{equation}
Each
$F_i: H_0^1(D)\to L^2([-1,1])$ is continuous, with the Fr\'{e}chet derivative and its adjoint given by \cite[Theorems 4.1 and 4.2]{Haltmeier07LK}
\begin{align}
    F_i^\prime(\xbs)h &= 2R_i(\xbs)R_i(h),\qquad\qquad \qquad\,\forall h \in H_0^1(D),\label{eqn:SC-derivative}\\
F_i^\prime(\xbs)^\ast g &= (I-\Delta)^{-1}(2R_i^\ast(R_i(\xbs)g)), \quad \forall g \in L^2([-1,1]),\label{eqn:SC-adjoint}
\end{align}
where $I$ is the identity operator, $\Delta$ denotes the Dirichlet Laplacian and $R_i^\ast: L^2([-1,1])\to L^2(D)$ denotes the adjoint of $R_i$ given by $R_i^\ast w(\boldsymbol{\xi}) := w(\langle \boldsymbol{\xi},\sigma_i\rangle)$.
Note, however, that the tangential cone condition has not been shown to hold for this problem.

Since $D\subset \mathbb{R}^2$ is bounded, by Sobolev embedding $H_0^1(D)\subset L^{r_\CX}(D)$ holds for $1<r_\CX<\infty$.
Thus, we can set $\CX=L^{r_\CX}(D)$, $r_\CX\in(1,\infty)$.
Moreover, since $L^2([-1,1])$ is continuously embedded into $L^{r_\CY}([-1,1])$ for $1<{r_\CY}\le2$, we may set $\CY_i = L^{r_\CY}([-1,1])$ for $r_\CY\in(1,2]$.
With these changes to the domain and range, operators $F_i:L^{r_\CX}(D)\to L^{r_\CY}([-1,1])$ with $r_\CX\in(1,\infty)$ and $r_\CY\in(1,2]$ are still Fr\'{e}chet differentiable, with the derivative and its adjoint given by \eqref{eqn:SC-derivative} and \eqref{eqn:SC-adjoint}. In the following, we omit the corresponding domain and use $L^r$ to denote both the $\CX$ and $\CY$ spaces.
Lebesgue spaces $L^{r}$ are known to be $r\vee2$-convex and $r\wedge 2$-smooth for $1<{r}<\infty$. The duality map $\CJ_p: L^r\to L^{r^\ast}$ can be computed by
\begin{align}\label{eqn:Lr_space_dualmap}
\CJ_p^{L^r}(\xbs) = \|\xbs\|_{L^r}^{p-r} |\xbs|^{r-1}\sign(\xbs).
\end{align}

\begin{figure}[hbt!]
    \centering
    \setlength{\tabcolsep}{0pt}
    \begin{tabular}{cccc}
    \includegraphics[height=.20\textwidth,trim={1.5cm .5cm .35cm 0cm},clip]{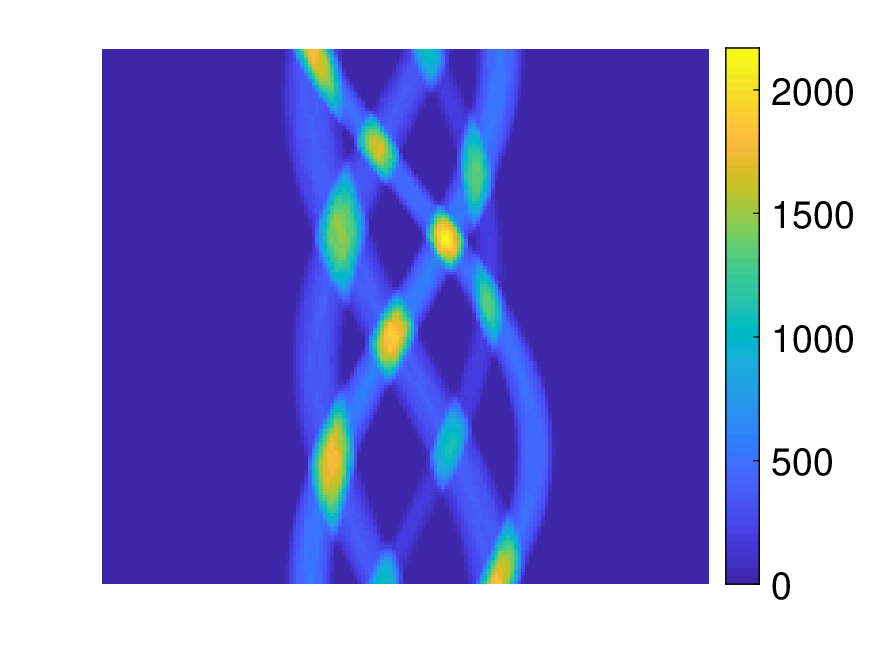}
   &\includegraphics[height=.20\textwidth,trim={1.5cm .5cm .35cm 0cm},clip]{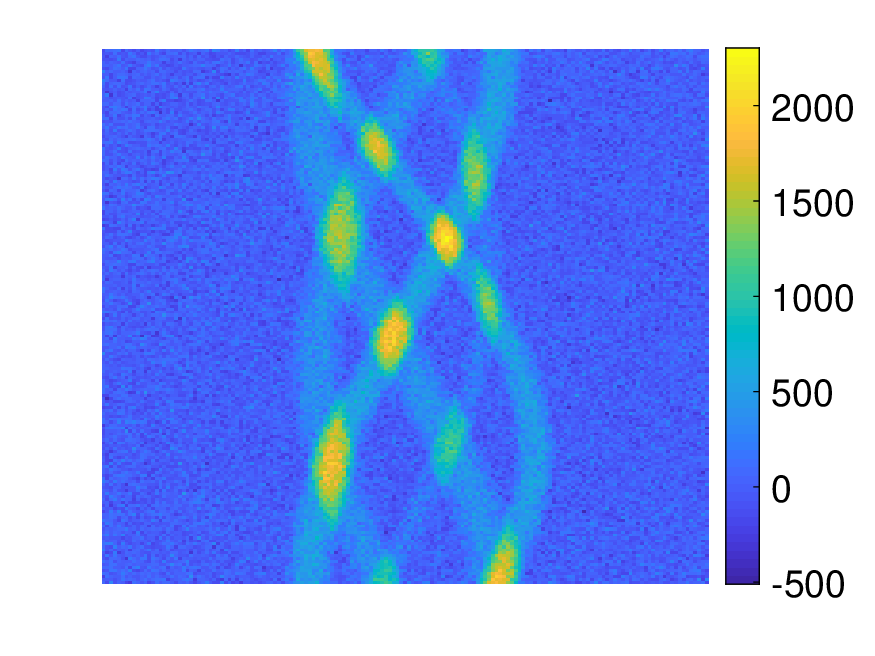}&
    \includegraphics[height=.20\textwidth,trim={1.5cm .5cm .35cm 0cm},clip]{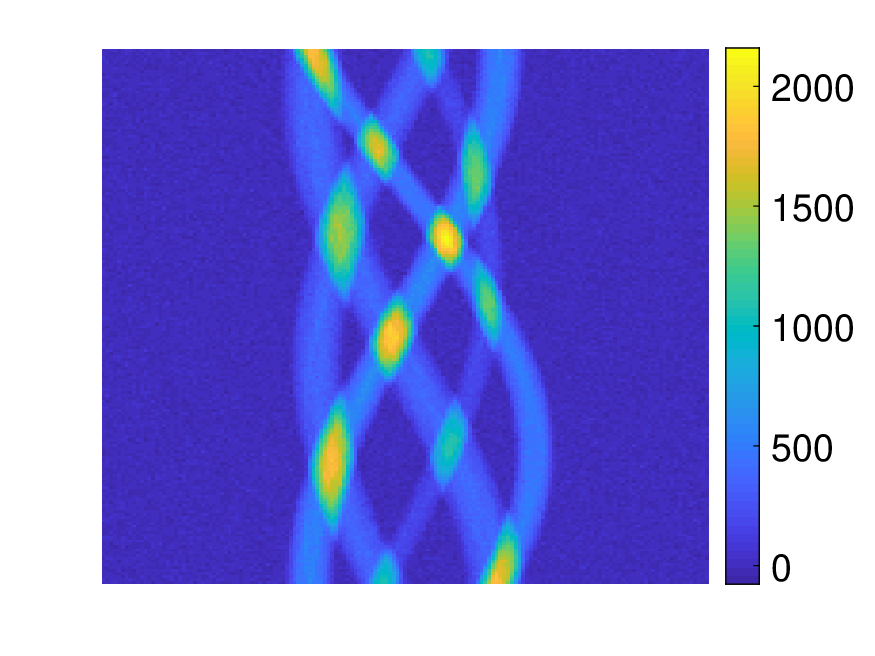}&
    \includegraphics[height=.20\textwidth,trim={1.5cm .5cm .35cm 0cm},clip]{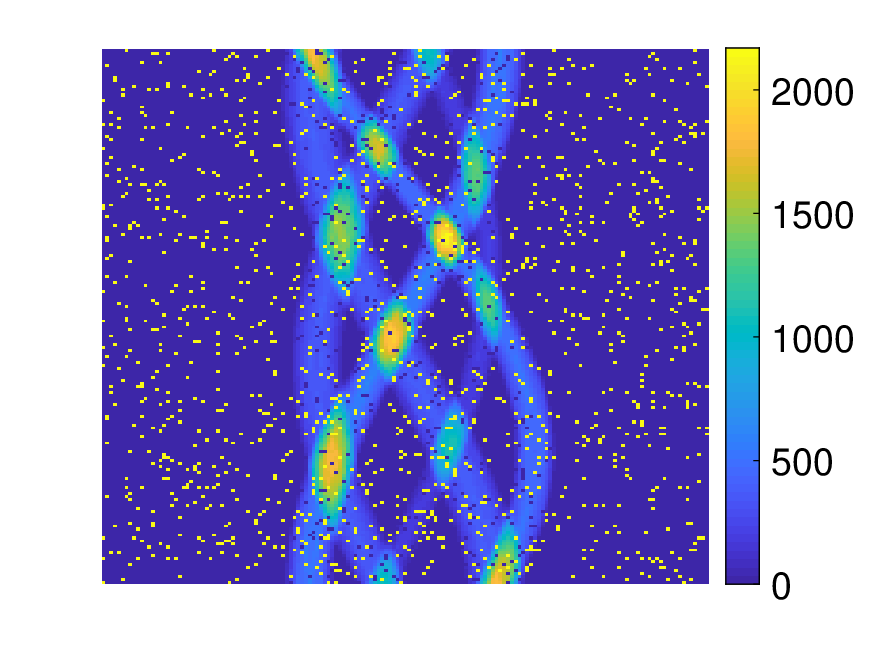}\\
    (a) exact data & (b) Gaussian $\epsilon = 0.05$ &
    (c) Gaussian $\epsilon =0.01$&
    (d) impulsive $\kappa = 0.1$
    \end{tabular}
    \caption{The noisy data  used in the recovery:
    {\rm(}b{\rm)}-{\rm(}c{\rm)} Gaussian noise;
    {\rm(}d{\rm)}: salt-and-pepper noise.}
    \label{fig:ST_data}
\end{figure}

To discretize the forward operator $F$, we partition the interval $[0,\pi)$ into $N=180$ equidistant angles  $\{0, \frac{\pi}{180},\ldots,\frac{179\pi}{180}\}$.
The sought-after solution $\xref$ is of dimension $110\times110$; see Fig. \ref{fig:ST_data}(a).
The forward operators $F_i$ are constructed with a batch size $b\geq 1$ by taking every $N/b$-th angle from the partition, starting with the $i$-th one.
The initial guess is set to $\xbs_0 = 0.01$, and the step-sizes $(\mu_k)_{k\in\mathbb{N}}$ follow the schedule $\mu_k =\mu_0 k^{-\alpha}$, with $\mu_0$ depending on the choice of $\CX$ and $\CY$, and $\alpha =0.2$, which satisfies the summability conditions $\sum_{k=1}^\infty \mu_k=\infty$ in Theorem \ref{thm:cone_convergence_qeqp}.

First, we study data with Gaussian noise. The noisy data $\ybs_i^\delta$, $i=1,\ldots,N$, are generated by
$\ybs_i^\delta = \ybs_i^\dag + \epsilon\|\ybs^\dag\|_{\ell^\infty}\xi_i$,
where $\|\cdot\|_{\ell^\infty}$ denotes the maximum norm of vectors, $\xi_i$s are i.i.d. and follow the standard Gaussian distribution, and $\epsilon>0$ is the noise level. Since the exact solution $\xbs^\dag$ is sparse, we  choose $1<r_\CX\le 2$. This leads to $p = 2$ under Assumption \ref{assn:X_and_Y}, and $q = 2$. With $\CY_i = L^2([-1,1])$, we compare three settings: $\CX=\CY=L^2$; $\CX=L^{1.5}$ and $\CY=L^{2}$; $\CX=L^{1.1}$ and $\CY=L^{2}$. SGD is run for 10\,000 epochs, where one epoch refers to $N/b$ iterations.

\begin{figure}[hbt!]
    \centering
    \begin{tabular}{cc}
    \includegraphics[height =.27\textwidth]{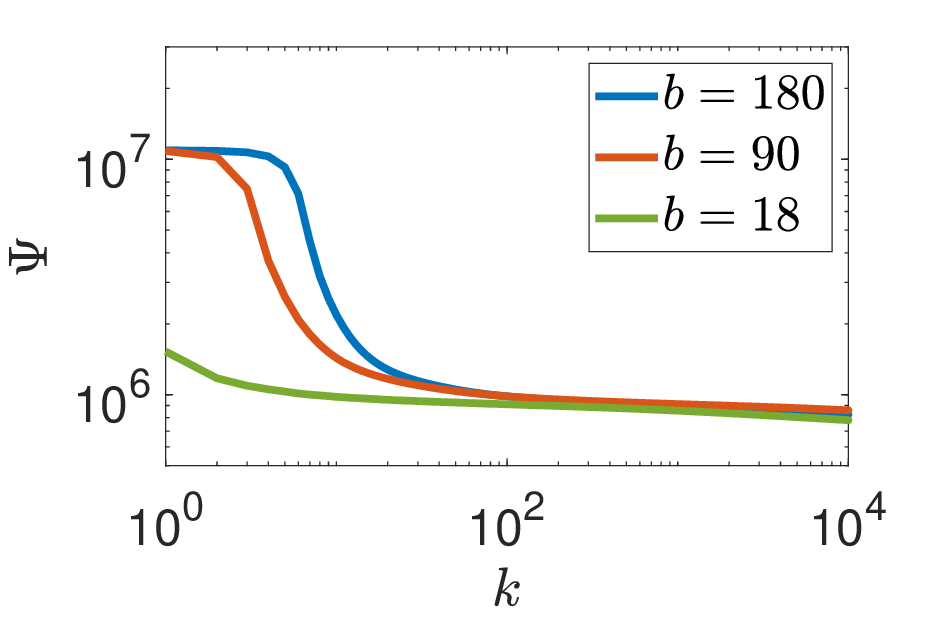}&
     \includegraphics[height =.27\textwidth]{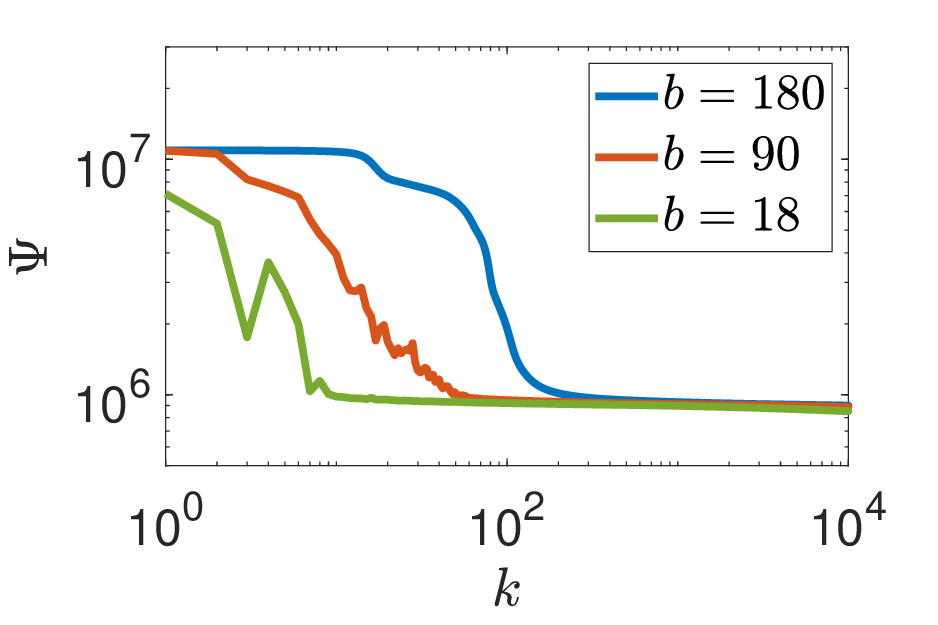}\\
    \includegraphics[height =.27\textwidth]{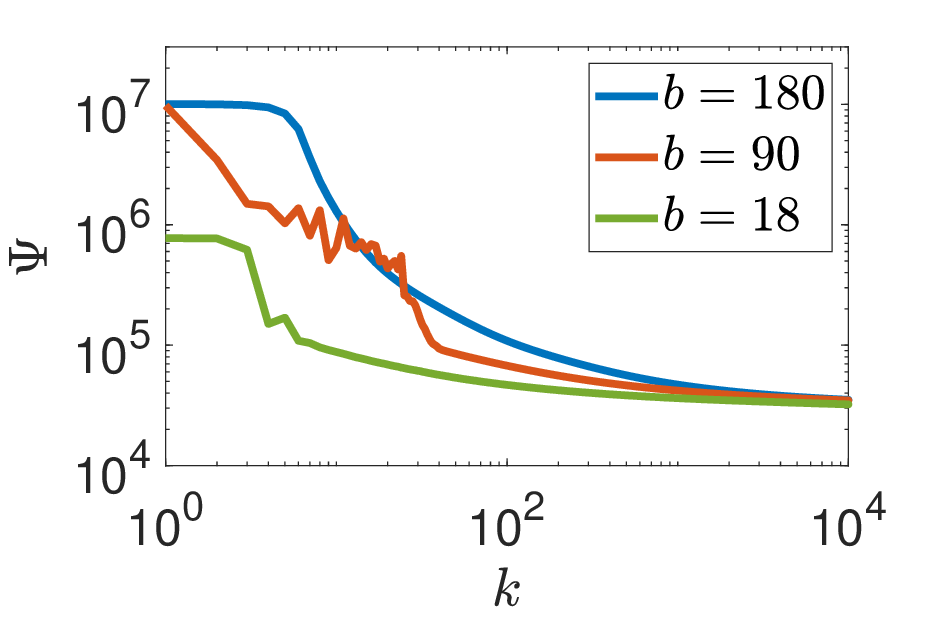}&
     \includegraphics[height =.27\textwidth]{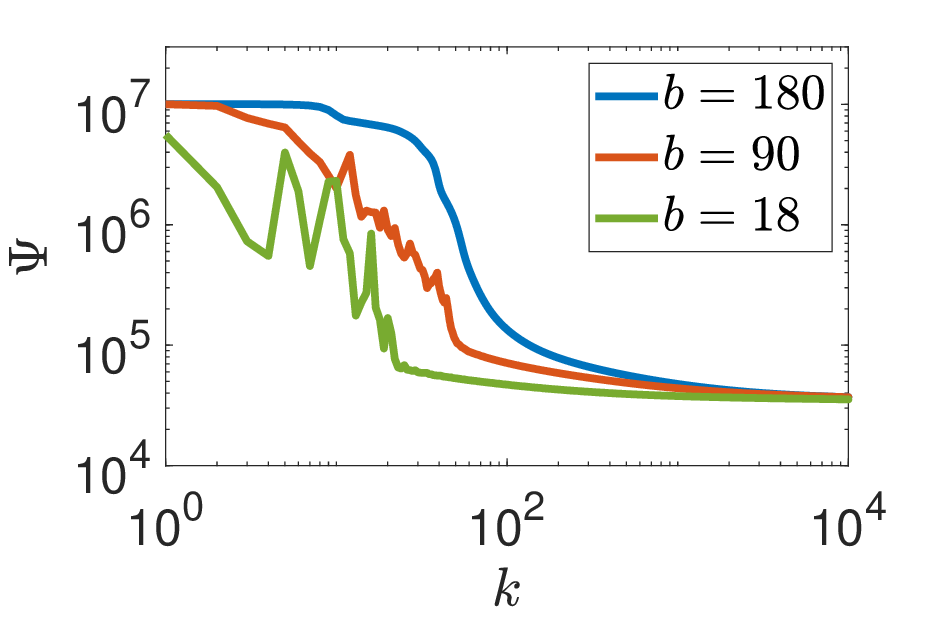}\\
    {{(a) $\CX=L^{1.5}$, $\CY=L^{2}$ }}&
    {{(b) $\CX=L^{1.1}$, $\CY=L^{2}$ }}
    \end{tabular}
    \caption{The objective $\Psi$ with respect to the batch size $b$. In the top row we use noise level $\epsilon = 5\times 10^{-2}$ and in the bottom row $\epsilon = 1\times 10^{-2}$.}
\label{fig:ST_Gaussian_res}
\end{figure}

In Fig. \ref{fig:ST_Gaussian_res}, we investigate the behavior of the objective $\Psi(\xbs_k)$, across iterates $\xbs_k$,  with respect to the batch size $b$.
We observe that a smaller batch size $b$ leads to a faster initial convergence, but is also accompanied by increased variance, as indicated by the oscillatory behavior.
Moreover, the results show that SGD with a smaller batch size $b$ is more effective than the Landweber method. Thus, we choose $b=18$ for the reconstruction task under Gaussian noise.

\begin{figure}[hbt!]
    \centering
    \begin{tabular}{cc}
    \includegraphics[height =.35\textwidth]{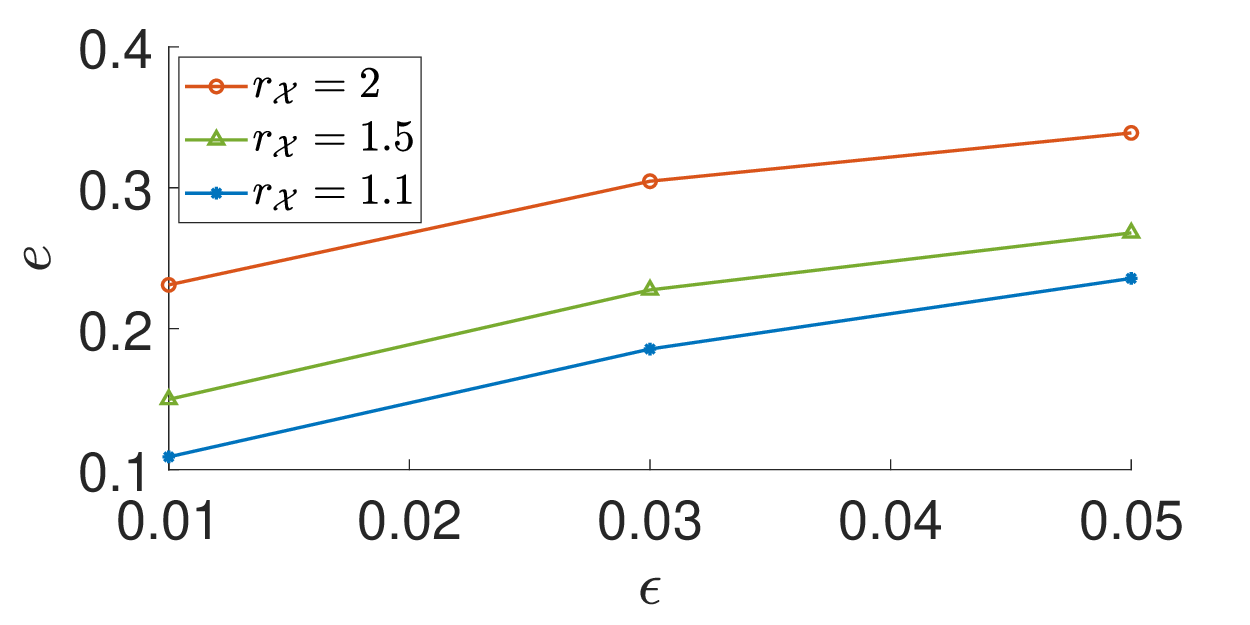}
    \end{tabular}
    \caption{The best reconstruction errors $e(\xbs)$ with respect to the noise level for three choices $\CX=L^{r_\CX}$.}
    \label{fig:ST_Gaussian_re}
\end{figure}

To measure the accuracy of the reconstruction, we adopt the relative $\ell^2$ error
$ e(\xbs)=\|\xref-\xbs\|_{\ell^2}/{\|\xref\|_{\ell^2}}$.
Fig. \ref{fig:ST_Gaussian_re} shows the error $e(\xbs)$ against the noise level $\epsilon$, with the error $e(\xbs)$ being minimal along the iteration trajectory.
The results are shown for three noise levels: $\epsilon = 5\times 10^{-2}$, $\epsilon = 3\times 10^{-2}$ and $\epsilon = 1\times 10^{-2}$. Note that the reconstruction errors are larger for $r_\CX = 2$ for all noise levels and decreases as $r_\CX$ tends to 1, i.e.,  the reconstruction quality improves as $r_\CX$ decreases.
These observations highlight the importance of selecting suitable $L^{r_\mathcal{X}}$ spaces.

\begin{figure}[hbt!]
    \centering
    \setlength{\tabcolsep}{0pt}
    \begin{tabular}{cccc}
        \includegraphics[height=.21\textwidth,trim={1.5cm .5cm 1cm 0cm},clip]{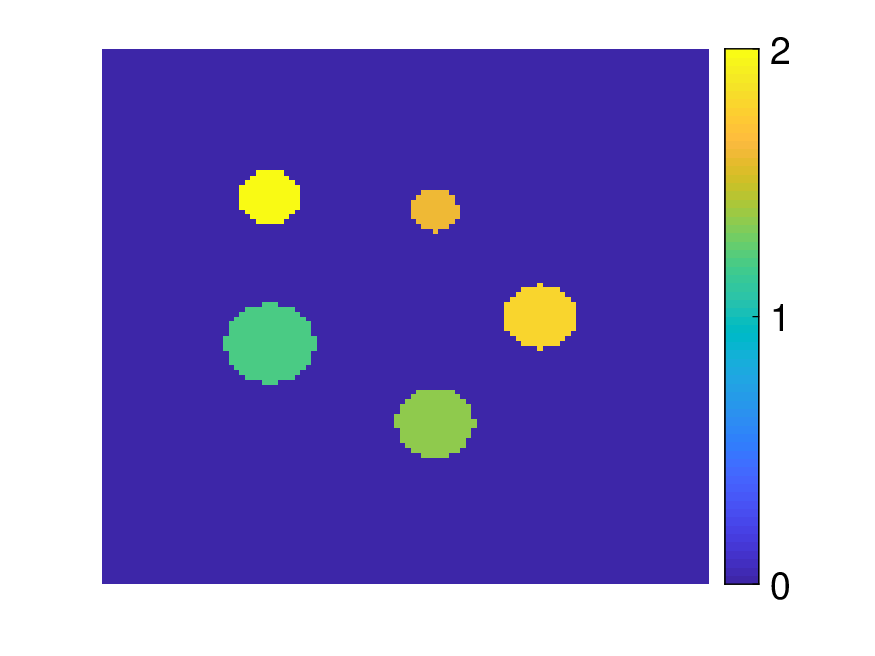}&
    \includegraphics[height=.21\textwidth,trim={1.5cm .5cm 1cm 0cm},clip]{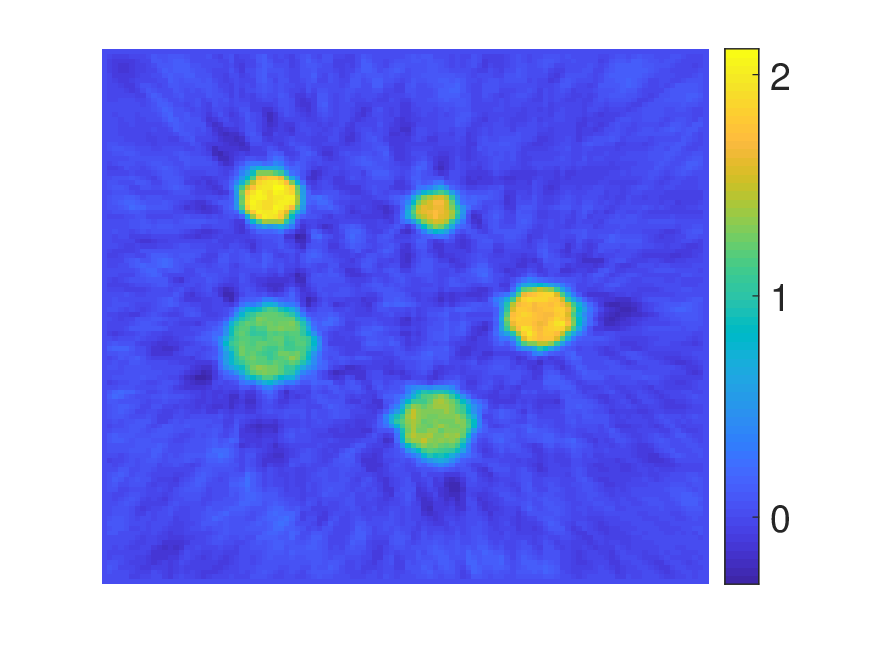}&
    \includegraphics[height=.21\textwidth,trim={1.5cm .5cm 1cm 0cm},clip]{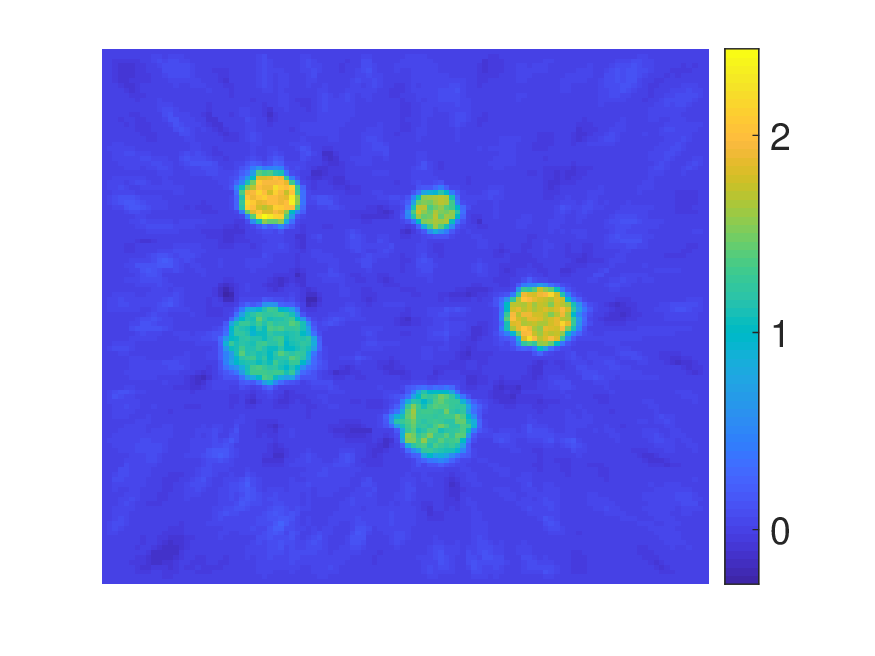}&
    \includegraphics[height=.21\textwidth,trim={1.5cm .5cm 1cm 0cm},clip]{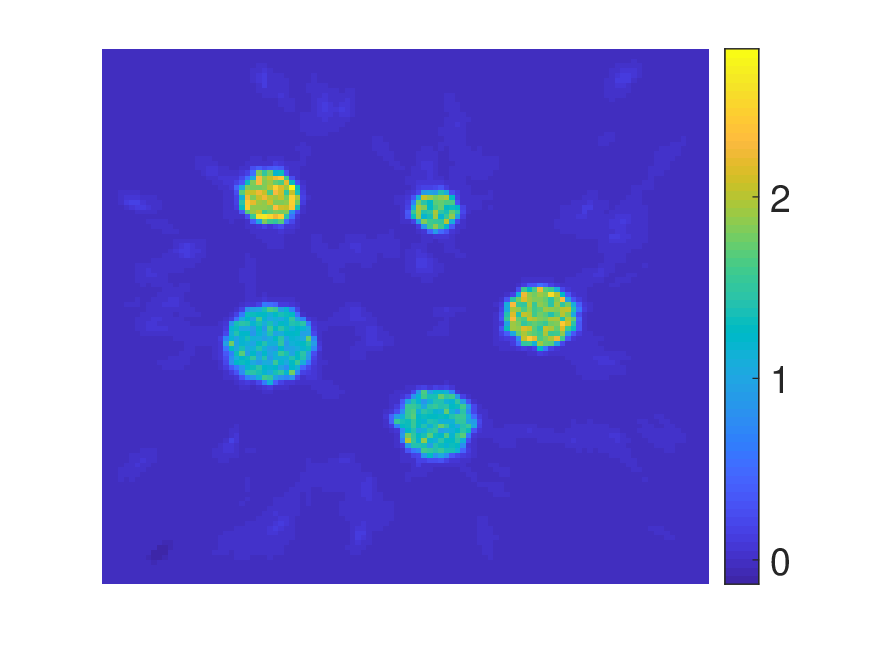}\\
    &{\scriptsize{(a) $\CX=\CY=L^2$ }}&
    {\scriptsize{(b) $\CX=L^{1.5}$, $\CY=L^{2}$ }}&
    {\scriptsize{(c) $\CX=L^{1.1}$, $\CY=L^{2}$ }}\\
    &\includegraphics[height=.21\textwidth,trim={1.5cm .5cm 1cm 0cm},clip]{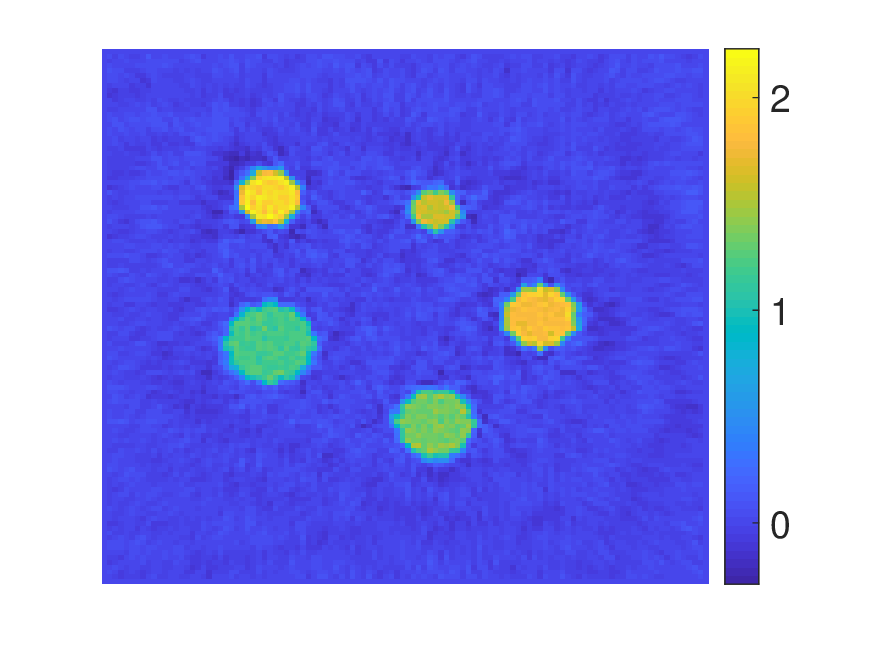} &
    \includegraphics[height=.21\textwidth,trim={1.5cm .5cm 1cm 0cm},clip]{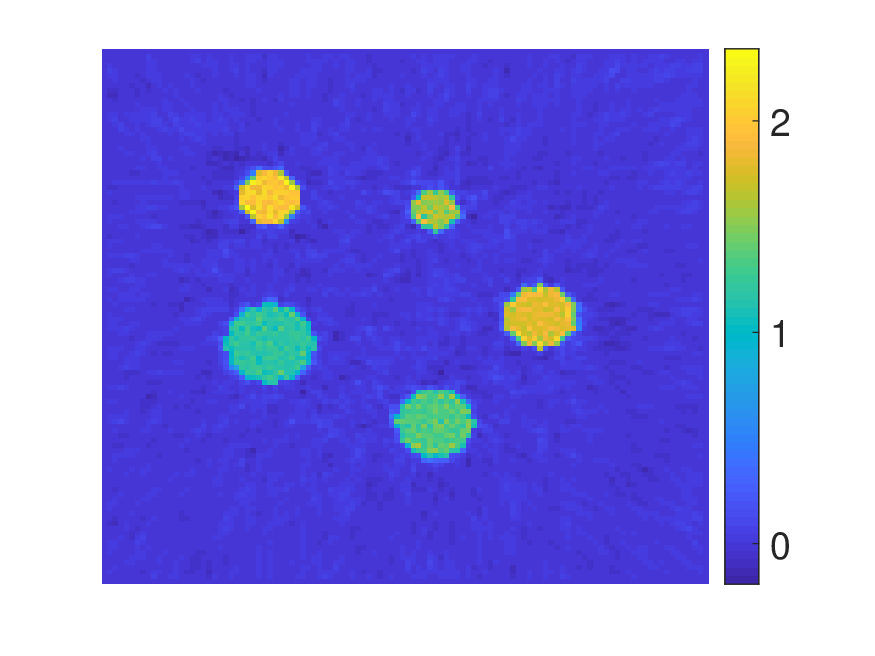}&
    \includegraphics[height=.21\textwidth,trim={1.5cm .5cm 1cm 0cm},clip]{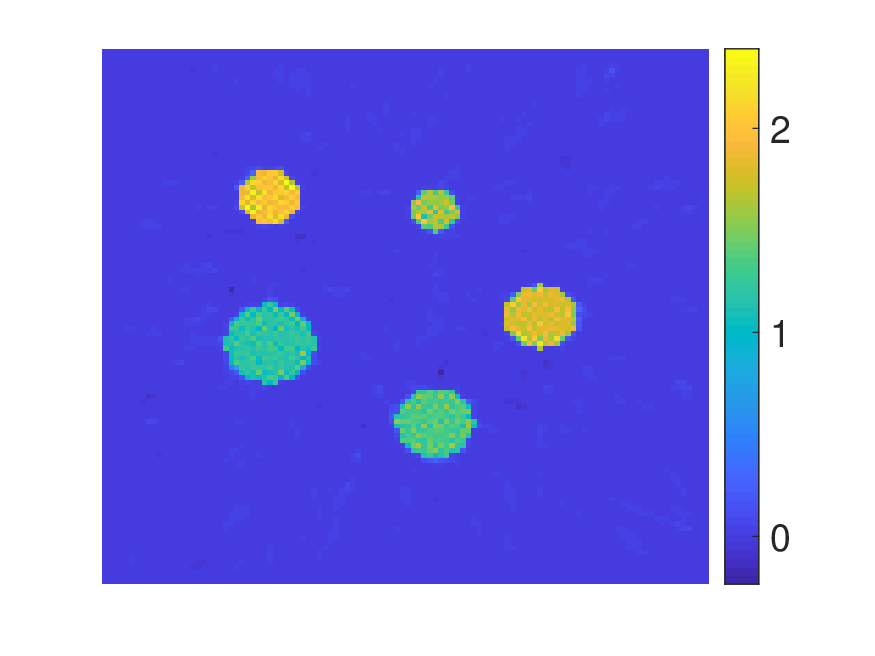}\\
    &{\scriptsize{(d) $\CX=\CY=L^2$ }}&
    {\scriptsize{(e) $\CX=L^{1.5}$, $\CY=L^{2}$ }}&
    {\scriptsize{(f) $\CX=L^{1.1}$, $\CY=L^{2}$ }}
    \end{tabular}
    \caption{The reconstruction results with batch size $b=18$.
    {\rm(}a{\rm)}-{\rm(}c{\rm)}: $\epsilon = 5\times 10^{-2}$.
    {\rm(}d{\rm)}-{\rm(}f{\rm)}: $\epsilon = 1\times 10^{-2}$.}
    \label{fig:ST_Gaussian}
\end{figure}

In Fig. \ref{fig:ST_Gaussian} we show the best reconstructions, corresponding to iterates achieving the smallest error $e(\xbs)$ over 10\,000 epochs, with noise levels $\epsilon=5\times 10^{-2}$ and $\epsilon=1\times 10^{-2}$. The standard Hilbert space setting is suitable for recovering smooth solutions, whereas other settings focus on sparse solutions. Note that $L^{r_\CX}$ spaces progressively enforce sparser solutions as the exponent $r_\CX$ tends to $1$.
Fig. \ref{fig:ST_Gaussian} shows that the setting with $\CX=L^{2}$ performs well at recovering the magnitude of nonzero entries, but is very poor at the recovery of the support of the solution, unlike Banach settings which perform much better.
The last column, $\CX=L^{1.1}$, achieves the best accuracy.
Lastly, as expected, the reconstruction quality improves with a decrease of the noise level $\epsilon$.

Next we investigate the data $\ybs_i^\delta$ with salt-and-pepper noise:
\begin{equation*}
    \ybs_i^\delta = \left\{\begin{array}{ccc}
        \ybs_i,      & {\rm with \,\, probability} \,\, 1-\kappa, \\
        \ybs_{\max}, & {\rm with \,\, probability} \,\, {\kappa}/{2}, \\
        \ybs_{\min}, & {\rm with \,\, probability} \,\, {\kappa}/{2},
    \end{array}\right.
\end{equation*}
where $\ybs_{\max}$ and $\ybs_{\min}$ are the maximum and minimum of the exact data $\ybs^\dag$, respectively, and $\kappa\in(0,1)$ denotes the percentage of corrupted data points. In the experiment, we fix $\kappa=0.1$, cf. Fig. \ref{fig:ST_data}(d), and choose $1<r_\CX\le 2$ (hence $q=p=2$) and set $\CY_i = L^{r_\CY}([-1,1])$ with $1<r_\CY\le 2$.
We compare four settings: (a) $\CX=\CY=L^{2}$; (b) $\CX=L^{2}$, $\CY=L^{1.1}$; (c) $\CX=L^{1.1}$, $\CY=L^{2}$; (d) $\CX=L^{1.1}$, $\CY=L^{1.1}$, and run $1\,{}000$ iterations.

\begin{figure}[hbt!]
    \centering
    \setlength{\tabcolsep}{0pt}
    \begin{tabular}{cccc}
    \includegraphics[height=.21\textwidth,trim={1.5cm .5cm 1cm 0cm},clip]{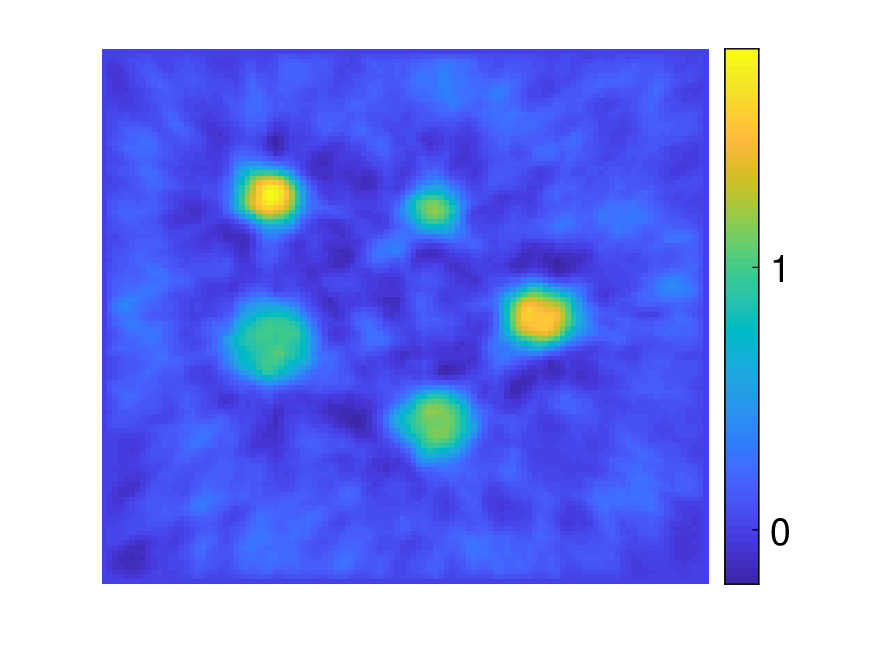}&
    \includegraphics[height=.21\textwidth,trim={1.5cm .5cm 1cm 0cm},clip]{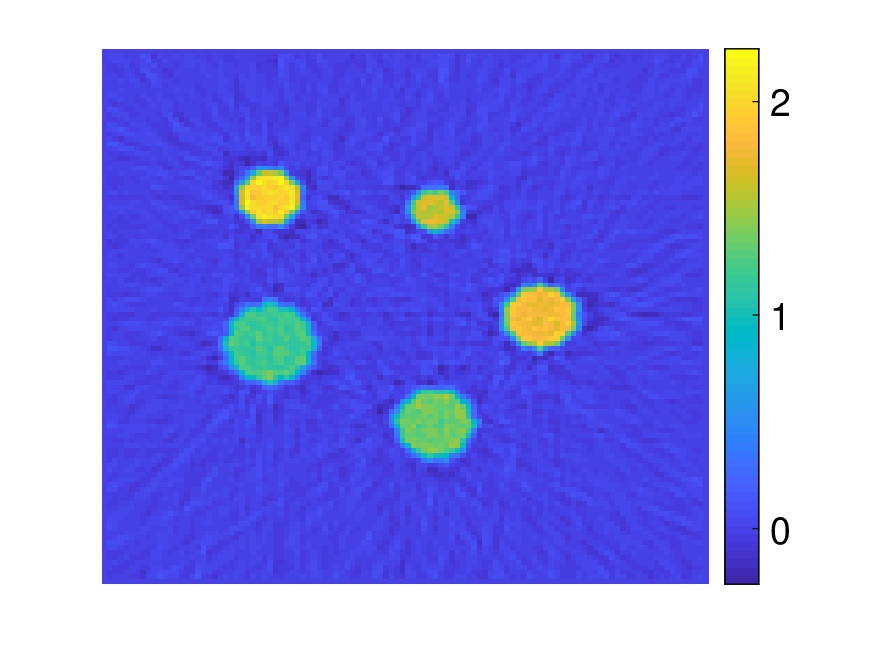}&
    \includegraphics[height=.21\textwidth,trim={1.5cm .5cm 1cm 0cm},clip]{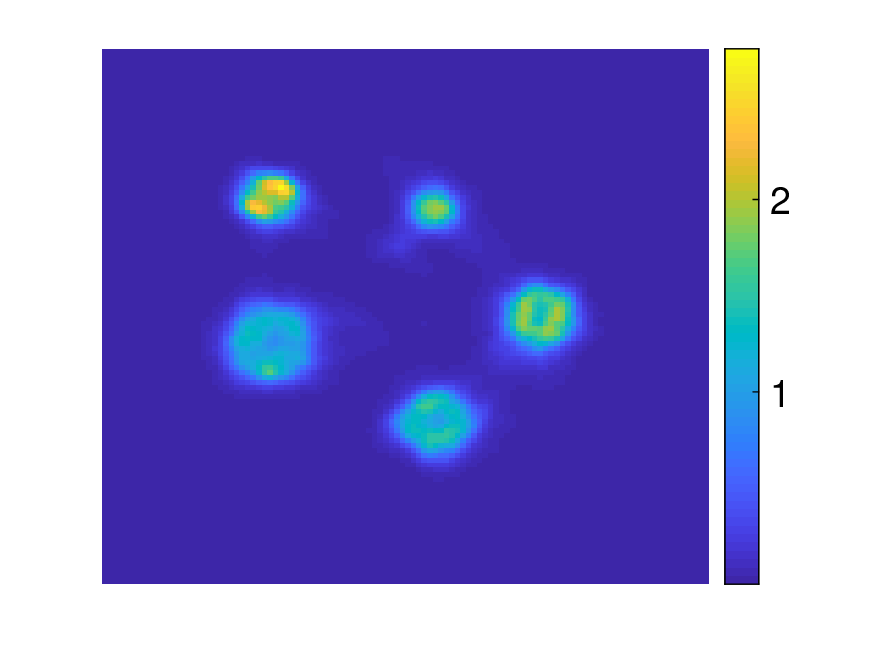} &
    \includegraphics[height=.21\textwidth,trim={1.5cm .5cm 1cm 0cm},clip]{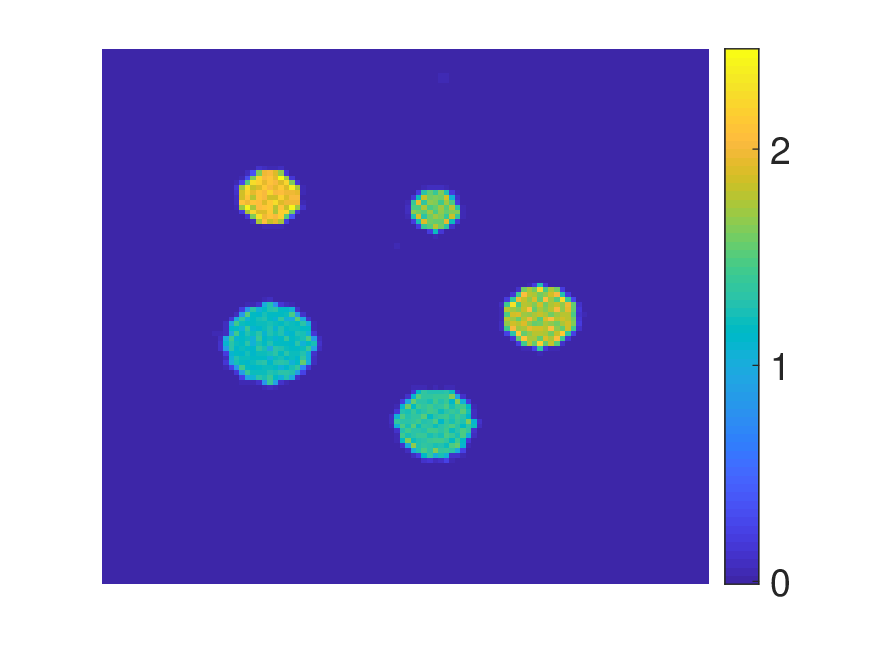}\\
    {\scriptsize{(a) $\CX=L^{2}$, $\CY=L^{2}$  }}&
    {\scriptsize{(b) $\CX=L^{2}$, $\CY=L^{1.1}$ }}& {\scriptsize{(c) $\CX=L^{1.1}$, $\CY=L^{2}$ }}&
    {\scriptsize{(d) $\CX=L^{1.1}$, $\CY=L^{1.1}$ }}\\
    {\scriptsize{$e(\xbs)$: $0.553$ }}&
    {\scriptsize{$e(\xbs)$: $0.230$ }}&
    {\scriptsize{$e(\xbs)$: $0.373$ }}&
    {\scriptsize{$e(\xbs)$: $0.123$ }}
    \end{tabular}
    \caption{The reconstructions of the phantom from the  data in Fig. \ref{fig:ST_data}(d), by SGD with $b=30$.}
    \label{fig:ST_SNP}
\end{figure}

We present the optimal reconstructions (along the iteration trajectory), and the corresponding error $e(\xbs)$, in Fig. \ref{fig:ST_SNP}.
Setting (a), aligned with the Hilbert space framework, falls short in achieving effective reconstructions. The choice $r_\CY = 1.1$ gives accurate recovery of the magnitudes of non-zero entries, but the presence of artifacts persists within the background. Setting (c) exhibits superior support recovery but fails at reconstructing the non-zero components. Setting (d) outperforms (a)-(c) and yields the most accurate results. Thus, $\CY = L^r$ with $r$ close to $1$ is performs the best for salt-and-pepper noise.

\subsection{Electrical impedance tomography (EIT)}\label{sec:eit}
Let $\Omega\subset\mathbb{R}^2$ be an open bounded domain with a Lipschitz boundary $\partial\Omega$. The EIT inverse problem aims to recover the conductivity $\sigma$ inside a domain $\Omega$ (e.g. a watertank) from voltage measurements on its boundary $\partial\Omega$. The conductivity is modeled as $\sigma \in L_{+}^\infty(\Omega):= \{v\in L^\infty(\Omega): v\ge C>0 \quad a.e. \,\,{\rm in} \,\, \Omega\} $.
The potential  $u$ satisfies the following elliptic PDE
\begin{equation}\label{eqn:EIT}
\left\{\begin{aligned}
-{\rm{div}}(\sigma\nabla u) &= 0,            \quad  {\rm{in}} \,\, \Omega,  \\
\sigma\partial_\nu u&= f,\quad \mbox{on }\partial\Omega,
\end{aligned}\right.
\end{equation}
where $\nu$ denotes the unit outward normal vector to the boundary $\partial\Omega$, and $f$ denotes the current applied on $\partial\Omega$ and is assumed to satisfy the conservation of charge $\int_{\partial\Omega}f{\rm d}s = 0$.
We normalize the solution $u$ by enforcing the grounding condition $\int_{\partial\Omega}u{\rm d}s = 0$, so that there exists a unique solution $u\in \widetilde{H}^1(\Omega):=\{v\in{H}^1(\Omega): \int_{\partial\Omega}v{\rm d}s = 0\}$.
In practice, one applies a finite number of currents $\{f_i\}_{i=1}^N$ and then records the resulting boundary potentials $\{g_i\}_{i=1}^N$. The forward operators $F_i$ are given by
\begin{equation}\label{eqn:EIT_forward_Kaczmarz}
    F_i : D(F)\subset L^\infty(\Omega)\to \tilde{H}^{1/2}(\partial\Omega), \quad \sigma\to g_i, \quad i= 1,\ldots, N.
\end{equation}
The Fr\'{e}chet differentiability of the operators $F_i$ in \eqref{eqn:EIT_forward_Kaczmarz} is well known \cite{LechleiterRieder2008,JinKhanMaass:2012}. For $\sigma\in D(F)$ and $\gamma\in L^\infty(\Omega)$, the Fr\'{e}chet derivative $F^\prime_i(\sigma)\gamma = w^i_\gamma|_{\partial\Omega}$, with $w_\gamma^i$ solving
\begin{equation*}
    \int_{\Omega}\sigma\nabla{w_\gamma^i}\cdot\nabla\varphi {\rm d}x= -\int_{\Omega}\gamma\nabla u_\sigma^i\cdot\nabla\varphi{\rm d}x, \quad \forall \varphi\in \tilde{H}^1(\Omega),
\end{equation*}
where $u_\sigma^i\in\tilde H^1(\Omega)$ is the solution of problem \eqref{eqn:EIT} with $f=f_i$. For $h\in \tilde{H}^{-1/2}(\partial\Omega)$, the adjoint of the derivative satisfies $F^\prime_i(\sigma)^\ast h = -\nabla u_\sigma^i\cdot \nabla\phi_h$, with the adjoint variable $\phi_h$, is obtained by solving
\begin{equation*}
    \int_{\Omega}\sigma\nabla{\phi_h}\cdot\nabla\varphi {\rm d}x= \int_{\partial\Omega}h\varphi{\rm d}s,\quad {\forall}  \varphi\in \tilde{H}^1(\Omega).
\end{equation*}
Since the space $L^\infty(\Omega)$ is not regular enough to be included in the analysis and $L^\infty(\Omega)\subset L^r(\Omega)$ for $r\in(1,\infty)$, we redefine the forward operator as
$F_i : D(F)\subset L^r(\Omega)\to \tilde{H}^{1/2}(\partial\Omega)$, $ \sigma\to g_i$.
However, the set $D(F)$ has no interior points in the $L^r$-topology. Thus, we restrict the sought-for conductivity space $\CX$ to a finite dimensional space $\mathcal{V}_r\subset L^\infty(\Omega)$, and redefine $F_i$'s as
$F_i : \mathcal{V}_{+}\subset \mathcal{V}_r\to L^{2}(\partial\Omega)$, $ \sigma\to g_i$,
with $\mathcal{V}_{+}:=  D(F)\cap\mathcal{V}$ and $\mathcal{V}_r:=(\mathcal{V},\|\cdot\|_{L^r})$.
By the norm equivalence in finite dimensional spaces, the differentiability of $F_i$ remains valid. The tangential cone condition for EIT in the general case remains open. Nonetheless, under suitable conditions on the measurements, the it is expected to hold for the operator $F$ in the admissible class of piecewise analytic conductivities; see  \cite[Sections 3 and 4]{Rieder2008_EIT} for in-depth discussions on the continuum  model and also the closely related complete electrode model for EIT. We refer interested readers to the work \cite{Kindermann:2022} for related results on the tangential cone condition for EIT using the tool of L\"{o}wner convexity.

Let $\Omega$ be the unit disc centered at the origin. We take $N$ currents, given as
$ f_n = \cos(\theta n)$, $n = 1,\ldots, N$ (with $N=20$ in the experiments), for $\theta\in[-\pi,\pi]$.
We employ the Galerkin finite element method (FEM) on a mesh with $M=4128$ triangle elements, which is also used to compute  Fr\'{e}chet derivative $F^\prime_i(\sigma)\gamma$ and its adjoint $F^\prime_i(\sigma)^\ast h$ and to reconstruct the conductivity $\sigma$.
With the piecewise constant basis functions $(v_j)_{j=1}^M$ on the mesh, the exact conductivity $\sigma$ is then given by
\begin{equation*}
    \sigma^{\dag}(x) = \bar\sigma(x)+\sum_{j=1}^M\alpha_jv_j(x), \quad
    \alpha_j = 3\chi_B,
\end{equation*}
where $\bar\sigma(x)\equiv 1$ denotes a known background.
We define the set $B$ as two open balls, with a radius $0.2$ and centers $(0,0.75)$ and $(0,-0.75)$, and $\chi_B$ denotes its characteristic function. The exact data $(g_i^\dag)_{i=1}^N$ are computed using the FEM with a finer mesh. The SGD is initialized with $\sigma_0=\bar\sigma\equiv 1$, and follows a step-size schedule $\mu_k =\mu_0 k^{-\alpha}$, with $\mu_0$ depending on the choice of $\CX$ and $\CY$, and $\alpha=0.05$. To measure the accuracy of the reconstruction, we employ the relative $L^2$ error $e (\sigma)= {\|\sigma-\sigma^\dagger\|_{L^2(\Omega)}}/{\|\sigma^\dagger\|_{L^2(\Omega)}}$.

\begin{figure}[hbt!]
    \centering
    \setlength{\tabcolsep}{0pt}
    \begin{tabular}{cccc}
    \includegraphics[height=.22\textwidth]{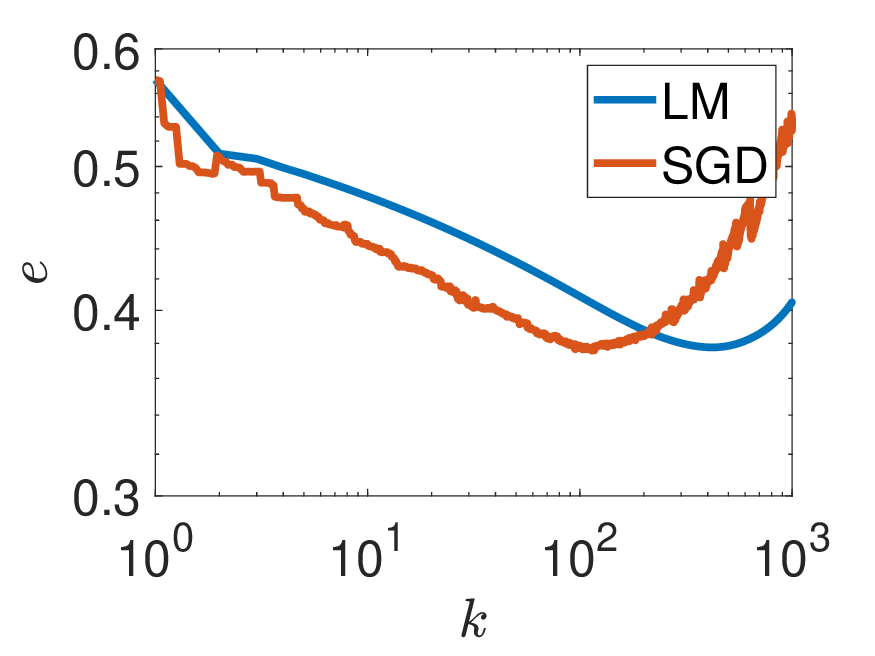}&
    \includegraphics[height=.22\textwidth]{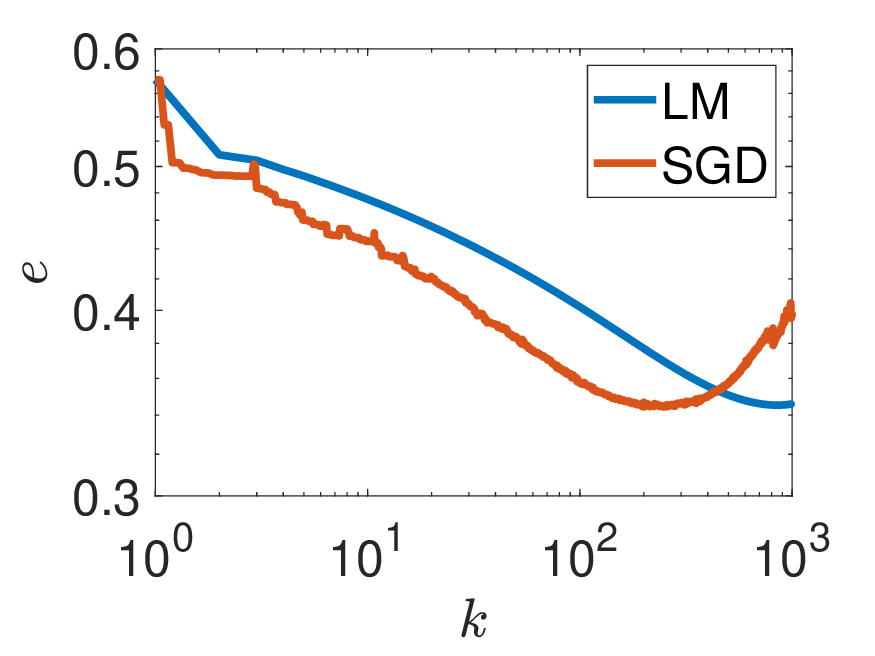}&
    \includegraphics[height=.22\textwidth]{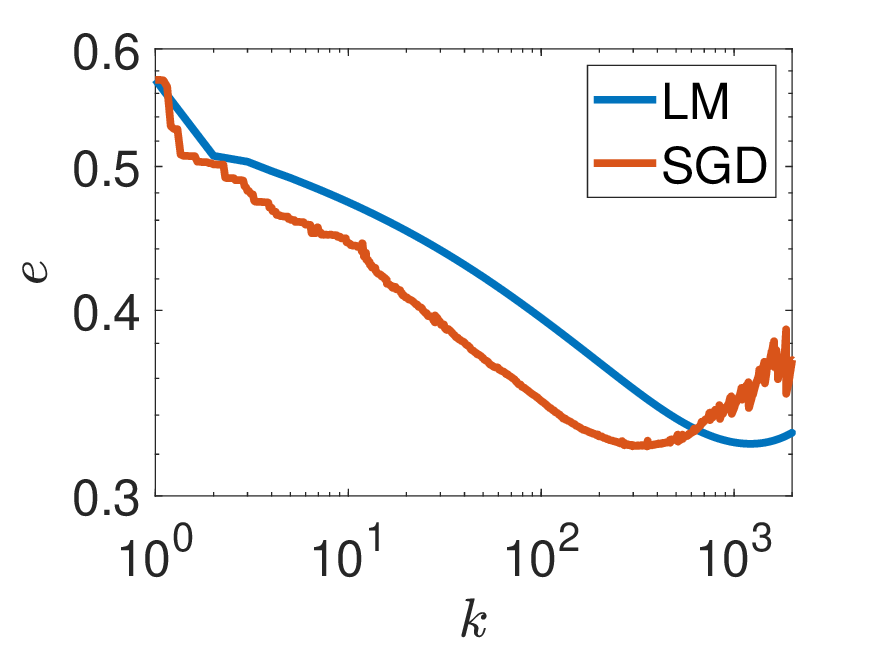}\\
    \includegraphics[height=.22\textwidth]{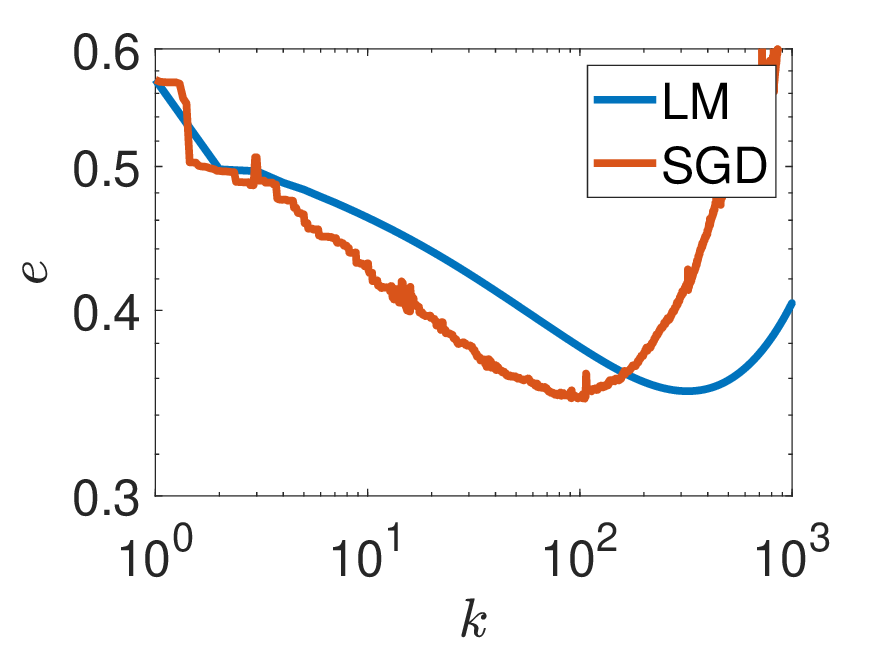}&
    \includegraphics[height=.22\textwidth]{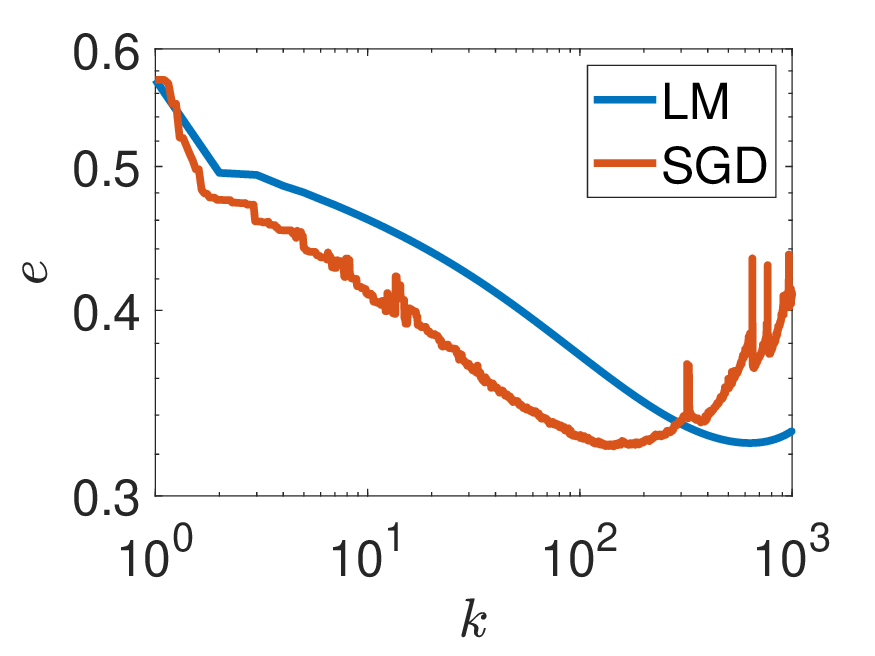}&
    \includegraphics[height=.22\textwidth]{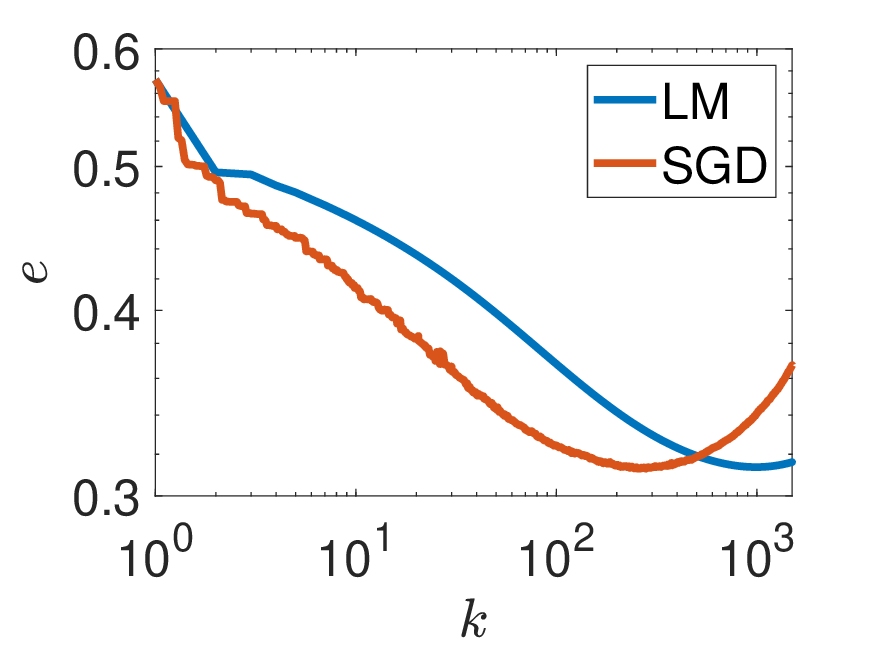}\\
    {\scriptsize{(a) $\epsilon = 2\times10^{-2}$ }}&
    {\scriptsize{(b) $\epsilon = 1\times10^{-2}$}}&
    {\scriptsize{(c) $\epsilon = 5\times10^{-3}$}}
    \end{tabular}
    \caption{The relative error $e(\sigma)$ versus the iteration number $k$ (in epoch) by LM and SGD in the case of noisy data. Top: $\CX=\CY=L^2$. Bottom: $\CX=L^{1.1}$ and $\CY=L^{2}$. }
    \label{fig:EIT_Gaussian_LM_vs_sgd}
\end{figure}

First we study the perturbed data $(g_i^\delta)_{i=1}^N$ formed by
$g_i^\delta = g_i^\dagger + \epsilon\max_{i=1,\ldots,N}\|g^\dagger_i\|_{L^\infty(\Omega)}\xi_i$,
where $\xi_i$s are i.i.d. and follow the standard Gaussian distribution and $\epsilon$ denotes a relative noise level. The function $\sigma^{\dagger}-\bar\sigma$ has a sparse representation. We take $1<r_\CX\le 2$, $r_\CY = 2$ and set $p=r_\CX$, $q=r_\CY$.
We investigate two Banach spaces settings. First with $\CX=\CY=L^2$, and second with $\CX=L^{1.1}$ and $\CY=L^{2}$.
In Fig. \ref{fig:EIT_Gaussian_LM_vs_sgd}, we compare SGD with the Landweber method (LM) for a maximum 1\,000 epochs, where one epoch refers to one Landweber iteration or $N$ SGD iterations. For both methods, the relative error $e(\sigma)$ first exhibits fast reduction, and then starts to increase, exhibiting a typical semi-convergence behavior. Moreover, the results in Fig. \ref{fig:EIT_Gaussian_LM_vs_sgd} show that SGD is more effective than the Landweber method, and the minimum errors are comparable for the two methods.
However, the inherent non-vanishing variance of stochastic gradients hinders the asymptotic convergence of SGD.

\begin{figure}[hbt!]
    \centering
    \setlength{\tabcolsep}{0pt}
    \begin{tabular}{cccc}
    &{\scriptsize{ $\epsilon = 2\times10^{-2}$ }}&
    {\scriptsize{ $\epsilon = 1\times10^{-2}$  }}&
    {\scriptsize{ $\epsilon = 5\times10^{-3}$  }} \\
    \includegraphics[height=.21\textwidth,trim={1.5cm .5cm 1cm 0cm},clip]{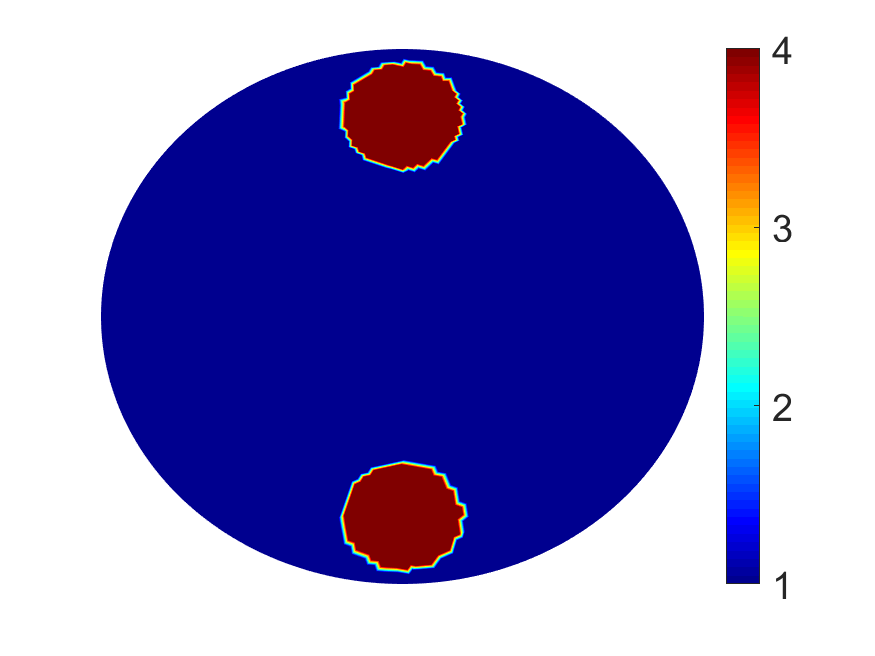}&
    \includegraphics[height=.21\textwidth,trim={1.5cm .5cm 1cm 0cm},clip]{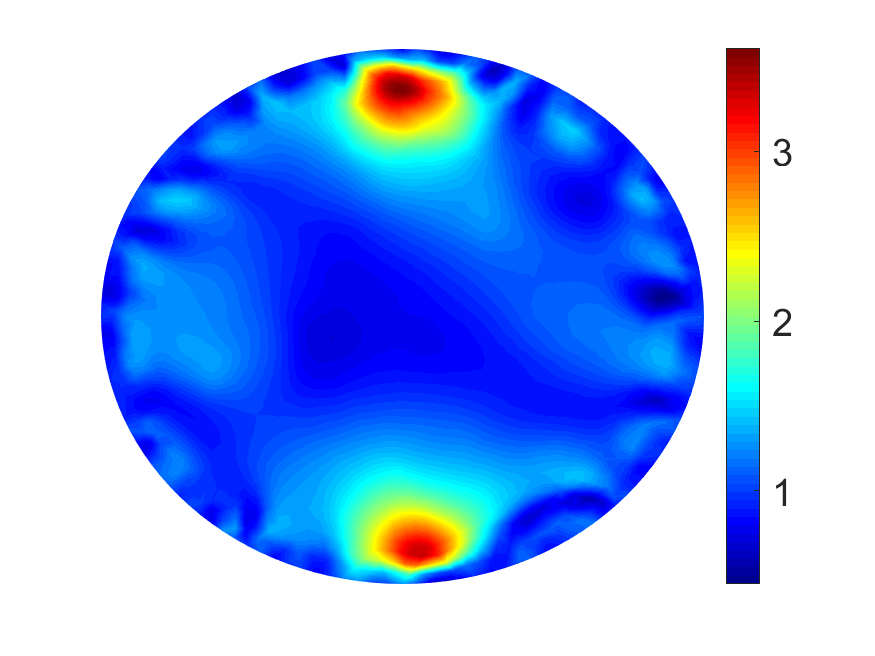}&
    \includegraphics[height=.21\textwidth,trim={1.5cm .5cm 1cm 0cm},clip]{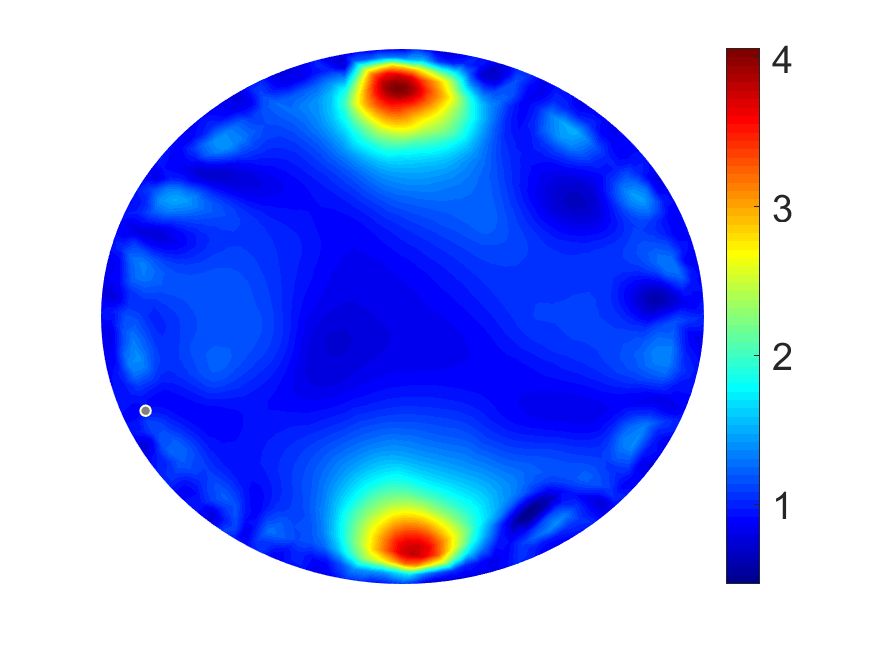}&
    \includegraphics[height=.21\textwidth,trim={1.5cm .5cm 1cm 0cm},clip]{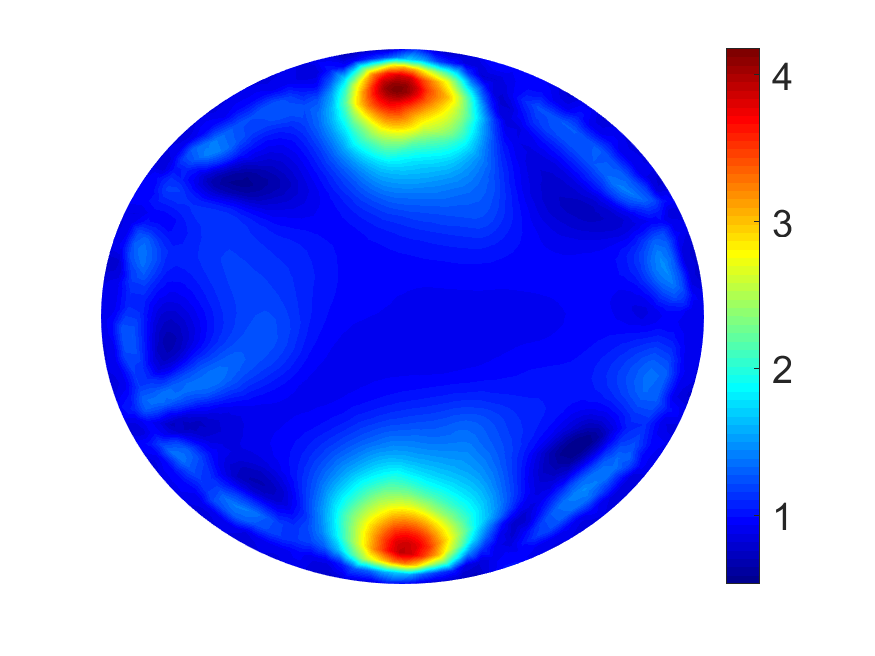}\\
    &{\scriptsize{(a) $\CX=L^{2}$, $\CY=L^{2}$ }}&
    {\scriptsize{(b) $\CX=L^{2}$, $\CY=L^{2}$ }}&
    {\scriptsize{(c) $\CX=L^{2}$, $\CY=L^{2}$ }}\\
    &{\scriptsize{ $e(\sigma) = 0.3758$ }}&
    {\scriptsize{ $e(\sigma) = 0.3444$ }}&
    {\scriptsize{ $e(\sigma) = 0.3239$ }}\\
    &\includegraphics[height=.21\textwidth,trim={1.5cm .5cm 1cm 0cm},clip]{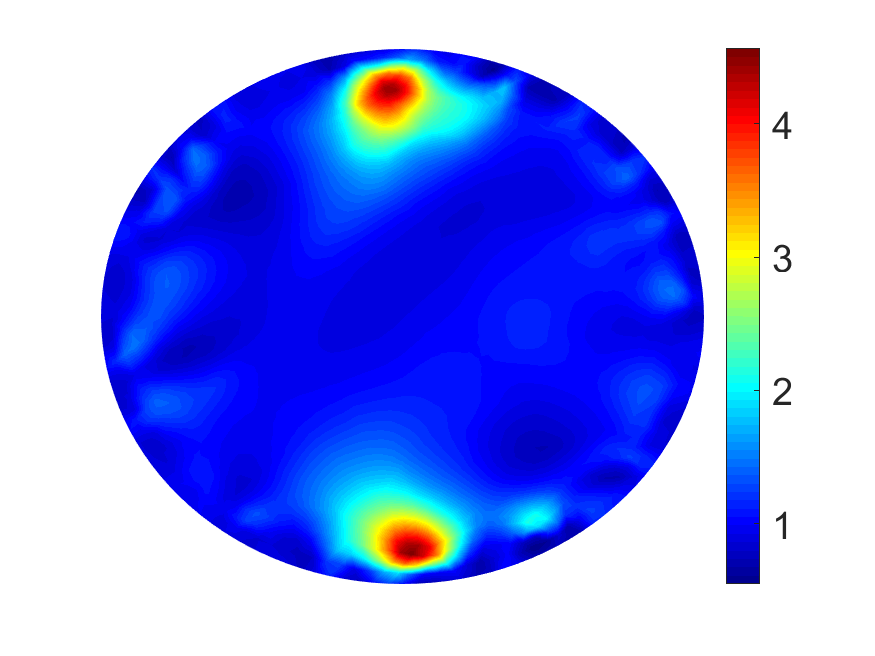}&
    \includegraphics[height=.21\textwidth,trim={1.5cm .5cm 1cm 0cm},clip]{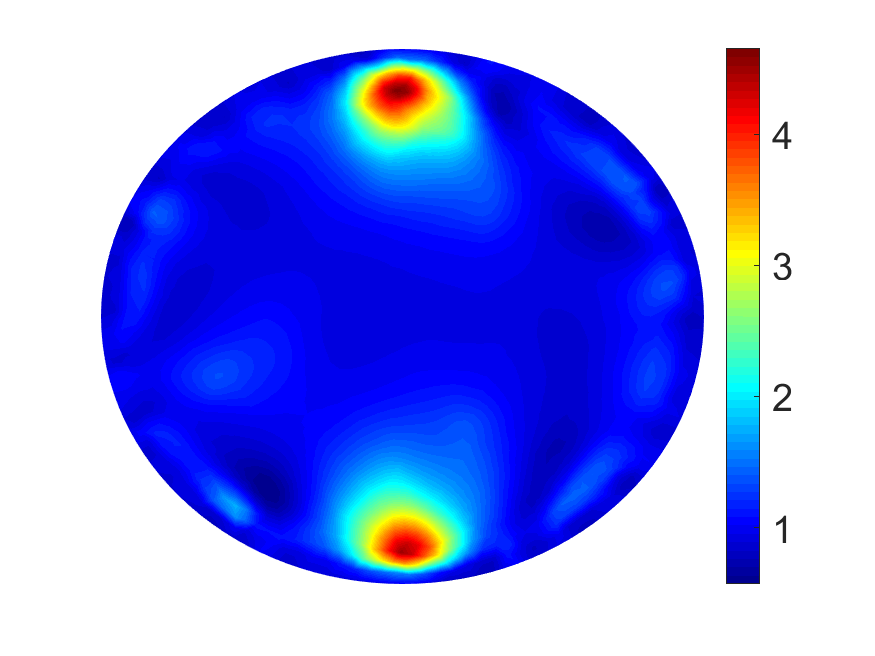}&
    \includegraphics[height=.21\textwidth,trim={1.5cm .5cm 1cm 0cm},clip]{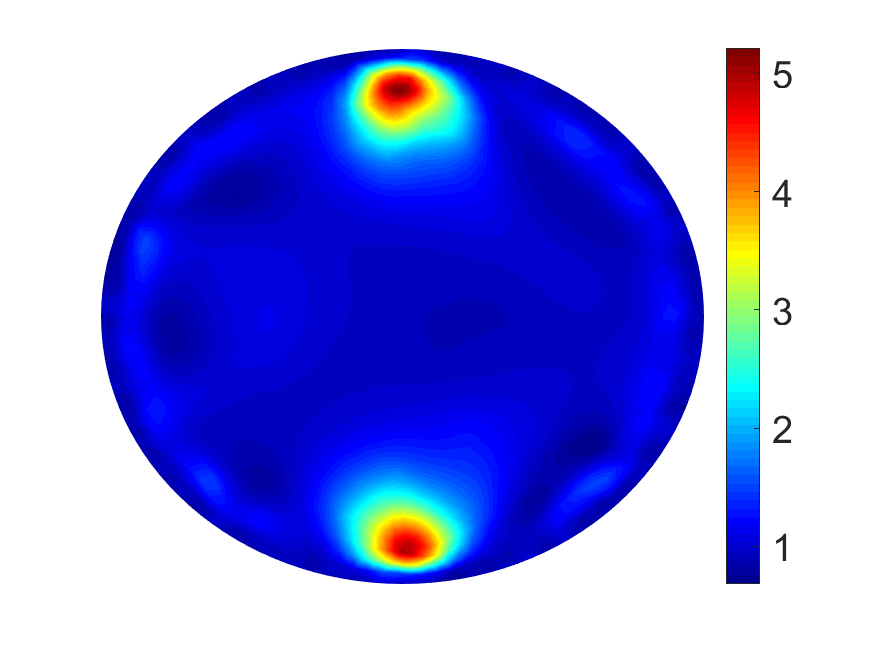}\\
    &{\scriptsize{(d) $\CX=L^{1.1}$, $\CY=L^{2}$ }}&
    {\scriptsize{(e) $\CX=L^{1.1}$, $\CY=L^{2}$ }}&
    {\scriptsize{(f) $\CX=L^{1.1}$, $\CY=L^{2}$ }}\\
    &{\scriptsize{ $e(\sigma) = 0.3488$ }}&
    {\scriptsize{ $e(\sigma) = 0.3239$ }}&
    {\scriptsize{ $e(\sigma) = 0.3128$ }}
    \end{tabular}
    \caption{The reconstructions by SGD for data with Gaussian noise.}
    \label{fig:EIT_noisy_data}
\end{figure}

To study the influence of Banach spaces, we show the optimal reconstructions (i.e., those with the smallest $e(\sigma)$) by SGD in Fig. \ref{fig:EIT_noisy_data}. The standard Hilbert space setting (row 1) yields smooth reconstructions. The Banach space setting (row 2) shows better reconstructions. This improvement agrees with intuition: the $L^{r_{\CX}}$ setting with the exponent $r_\CX$ close to 1 is more effective at recovering sparse solutions \cite{DDD04}.

\begin{figure}[hbt!]
    \centering
    \setlength{\tabcolsep}{0pt}
    \begin{tabular}{cccc}
     \includegraphics[height=.21\textwidth,trim={1.5cm .5cm 1cm 0cm},clip]{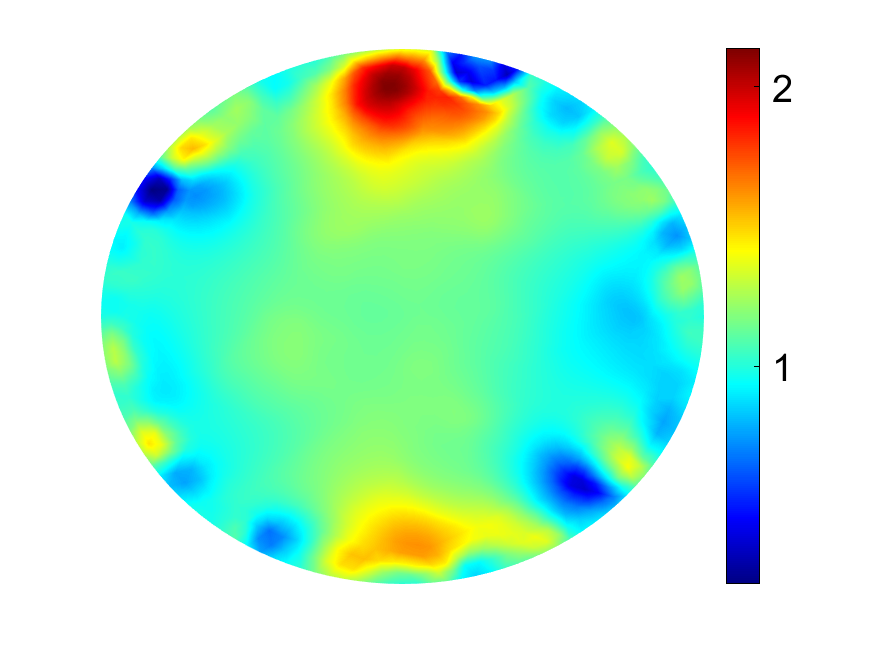}
    &\includegraphics[height=.21\textwidth,trim={1.5cm .5cm 1cm 0cm},clip]{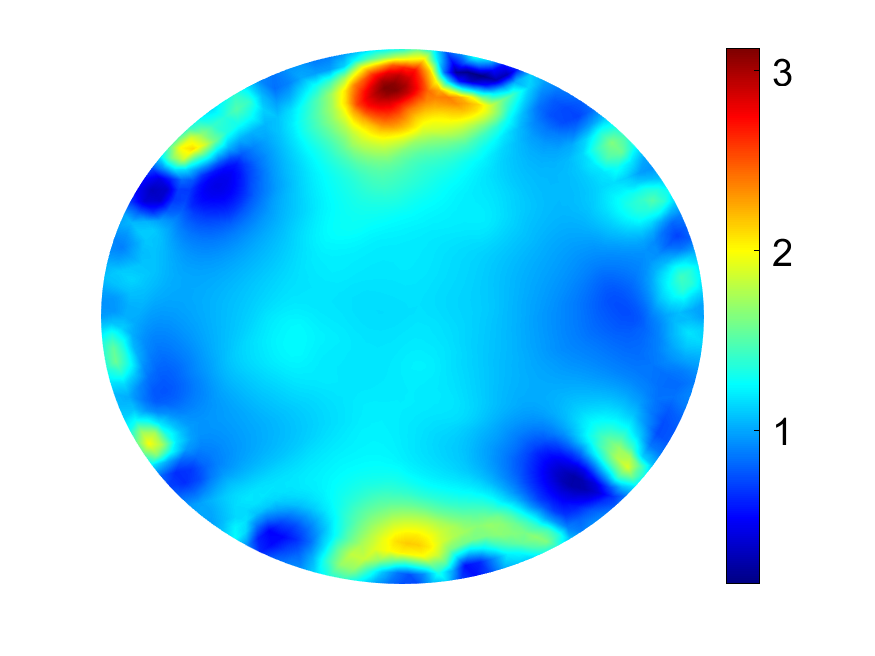}
    &\includegraphics[height=.21\textwidth,trim={1.5cm .5cm 1cm 0cm},clip]{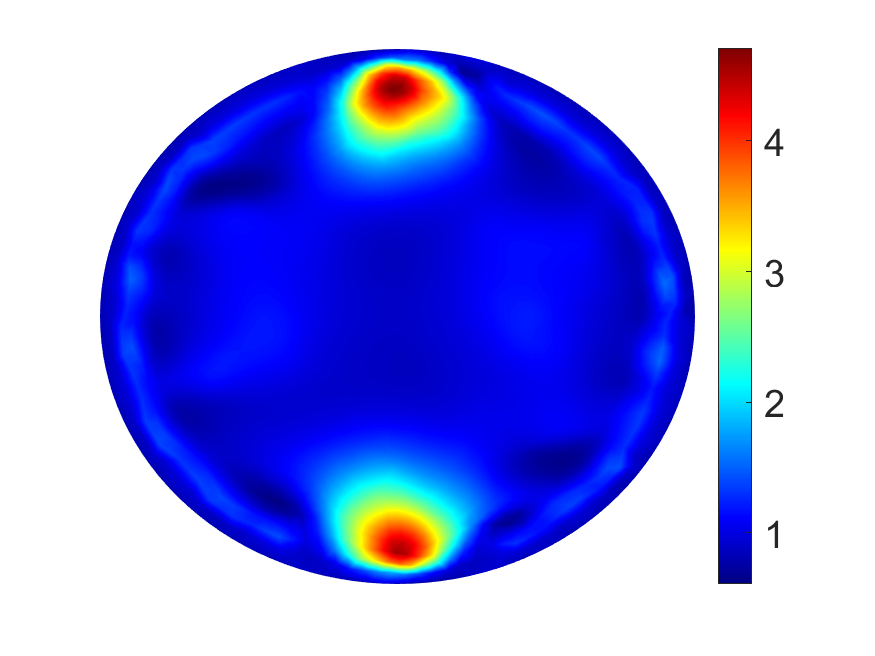}
    &\includegraphics[height=.21\textwidth,trim={1.5cm .5cm 1cm 0cm},clip]{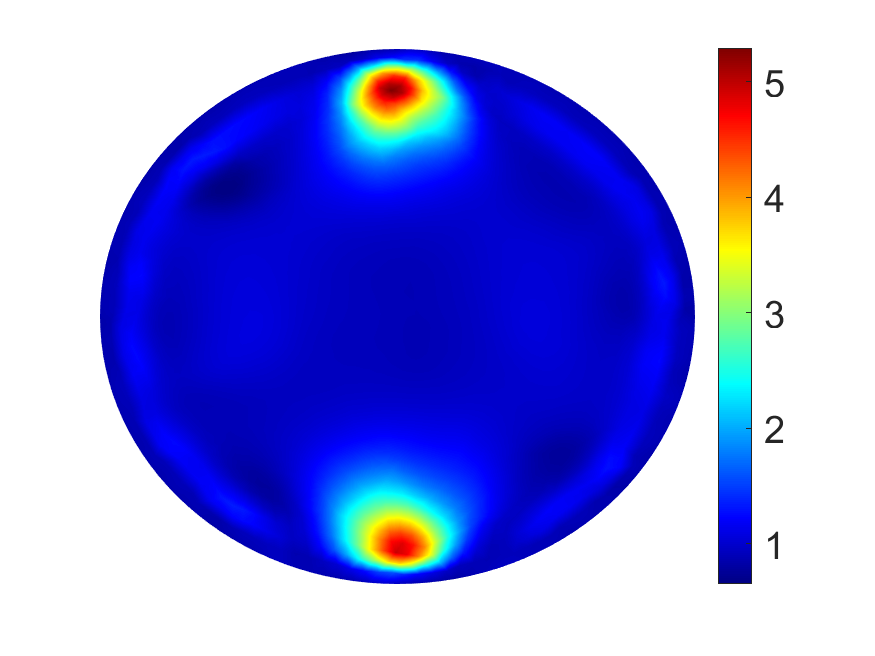}
    \\
    {\scriptsize{(a) $\CX=\CY=L^2$ }}&
    {\scriptsize{(b) $\CX=L^{1.1}$, $\CY=L^{2}$ }}&
    {\scriptsize{(c) $\CX=L^{2}$, $\CY=L^{1.1}$}}&
    {\scriptsize{(d) $\CX=L^{1.1}$, $\CY=L^{1.1}$}}\\
    {\scriptsize{ $e(\sigma) = 0.4829$ }}&
    {\scriptsize{ $e(\sigma) = 0.4497$ }}&
    {\scriptsize{ $e(\sigma) = 0.3050$ }}&
    {\scriptsize{ $e(\sigma) = 0.3086$ }}
    \end{tabular}
    \caption{The reconstruction results by SGD for data with impulsive noise.}
    \label{fig:EIT_impulsive}
\end{figure}

Next we study the data with impulsive noise:
\begin{equation*}
    g_i^\delta = \left\{
    \begin{array}{cc}
        g_i^\dagger, &  {\rm with \,\, probability} \,\, 1-\kappa,\\
        g_i^\dagger + \epsilon\max_{i=1,\ldots,N}\|g_i^\dagger\|_{L^\infty(\Omega)}\xi_i, & {\rm with \,\, probability} \,\, \kappa,
    \end{array} \right.
\end{equation*}
where $\xi_i$s follow the standard Gaussian distribution, $\kappa\in(0,1)$ is the percentage of corrupted data points and $\epsilon$ denotes the magnitude of noise. Below we fix $\kappa = 0.1$ and $\epsilon=0.4$. Since the $L^r$ fitting with $r$ close to $1$ is suitable for impulsive noise, we choose $1<r_\CX\le 2$ and $1<r_\CY\le 2$, and also $p=r_\CX$ and $q=r_\CY$.
Fig. \ref{fig:EIT_impulsive} presents the reconstructions,  obtained by running the SGD iteration and then selecting the iterate with the smallest $e(\sigma)$. The choice $\CY = L^{1.1}$ leads to a superior quality than that in Figs. \ref{fig:EIT_impulsive}(a) and (b). Thus, the $L^r$ space with $r\approx1$ is more appropriate to deal with the impulsive noise than the $L^2$ space. The comparison with Figs. \ref{fig:EIT_impulsive}(a) and (b) indicates that a small value for $r_\CX$ performs better, since it gives less variation in the background.

\begin{remark}
By Assumption \ref{assn:X_and_Y}, the parameter $p$ of the duality map $\CJ_p^{L^{r_\CX}}$ in \eqref{eqn:sgd} must be set to $2$ for $1<r_\CX\le2$, and the analysis requires $q = 2$, for which $\CJ_p^{L^{r_\CX}}$ and $\CJ_q^{L^{r_\CY}}$ can be computed using \eqref{eqn:Lr_space_dualmap}. However, the reconstructions with this choice are not as good as the ones using the duality map $\CJ_{r}$. This was also observed in \cite[Section 5.1.2]{Margotti2015}. Thus, we have chosen $p=r_\CX$ and $q=r_\CY$, and  $\CJ_{r}$ in the space $L^r$ is given by $\CJ_r^{L^r}(\xbs) =  |\xbs|^{r-1}\sign(\xbs)$. Note that, in principle, $p$ is not allowed when $p < 2$, since $\CX$ is assumed to be $p$-convex, for which there is still no convergence guarantee. Nevertheless, the duality map $\CJ_{r}$ is a natural choice for the space $L^r$ in practice. Developing a complete theory for the setting is an interesting future research problem.
\end{remark}

\appendix
\section{Auxiliary results}
\begin{theorem}[{Robbins-Siegmund theorem on the convergence of almost super-martingales, \cite[Lemma 11]{P87}}] \label{thm:sub_martingale_convergence}
Consider a filtration $(\CF_k)_{k\in\bbN}$ and four non-negative, $(\CF_k)_{k\in\bbN}$ adapted processes $(\alpha_k)_{k\in\bbN}$, $(\beta_k)_{k\in\bbN}$, $(\gamma_k)_{k\in\bbN}$, and $(\delta_k)_{k\in\bbN}$, such that $\sum_{k=1}^\infty \beta_k<\infty$ and $\sum_{k=1}^\infty \gamma_k<\infty$ with probability $1$, and for all $k$, we have
\[\bbE_k[\alpha_{k+1}]\le(1+\beta_k)\alpha_k+\gamma_k-\delta_k.\]
Then the sequence $(\alpha_k)_{k\in\bbN}$ converges a.s. to a random variable $\alpha_\infty$, and $\sum_{k=1}^\infty\delta_k <\infty$ a.s.
\end{theorem}

\begin{lemma}[Polyak inequality, {\cite[Lemma 6]{P87}}]\label{lem:polyak_series}
Let $(\delta_n)_n$ be a sequence of non-negative scalars, $(\mu_n)_n$ a sequence of positive scalars, and $\alpha>0$.
If
$\delta_{n+1}\leq \delta_n-\mu_{n+1}\delta_n^{1+\alpha}$ for any $n=0,\ldots,N$,
then
\[
\delta_N\leq {\delta_0}{\Big(1+\alpha\delta_0^\alpha\sum_{n=1}^{N}\mu_n \Big)^{-1/\alpha}}.
\]
\end{lemma}

\section{The proof of Theorem \ref{thm:cone_convergence_qeqp}}\label{app:cone_convergence_qeqp}
\begin{proof}
By taking the conditional expectation on \eqref{eqn:bregman_descent} with respect to $i_k$, we get
\begin{align}\label{eqn:cond_expect_bregman}
\bbE_k[\Delta_{k+1}]
\leq\Delta_k -p \Big(1-\gamma -\Lmax^{p^\ast}\frac{G_{p^\ast}}{p^\ast}\mu_{k+1}^{p^\ast-1}\Big)\mu_{k+1}\Psi(\xbs_k).
\end{align}
Using the assumed bound on $\mu_k$ in Lemma \ref{lem:bregman_monotonicity_qeqp} and taking the full expectation yield
\begin{align*}
\bbE[\Delta_{k+1}]\leq\bbE[\Delta_k] -p C\mu_{k+1}\bbE[\Psi(\xbs_k)].
\end{align*}
Thus, $(\Delta_k)_{k\in\bbN}$ is a non-negative super-martingale (cf. Theorem \ref{thm:sub_martingale_convergence}) and consequently it a.s. converges to some non-negative variable $\Delta_\infty$, such that $\sum_{k=0}^\infty \mu_{k+1}\Psi(\xbs_k)<\infty$ on an a.s. set $\Omega$.
Here, $\Omega$ is the measurement set on which $(\Delta_k)_{k\in\bbN}$ is convergent, $\sum_{k=0}^\infty \mu_{k+1}\Psi(\xbs_k)<\infty$, and $\bbP(\Omega)=1$.
Moreover, by \cite[Lemma 5.31]{BauschkeCombettes:2011}, it follows that $\lim_{k\rightarrow\infty} \mu_k\bbE[\Psi(\xbs_k)]=0$.
Thus, if there exists a $k_0$ and a constant $M$ such that $\mu_k\geq M$ for all $k\geq k_0$, it follows that $\lim_{k\rightarrow\infty} \Psi(\xbs_k)=0$ a.s. and $\lim_{k\rightarrow\infty} \bbE[\Psi(\xbs_k)]=0$.
Otherwise, since $\sum_{k=0}^\infty\mu_k=\infty$, we have $\liminf_k\Psi(\xbs_k)=0$ for every event $\omega\in\Omega$, which can be proven by constructing a contradiction.
Consider an event $\omega$ on which this is not the case, i.e., $\liminf_k\Psi(\xbs_k)>0$. Then there exists $\epsilon>0$ and $k_\epsilon\in\bbN$ such that for all $k\ge k_\epsilon$, $\Psi(\xbs_k)\ge\epsilon$, giving $\sum_{k\ge k_\epsilon}\mu_{k+1}\Psi(\xbs_k)\ge\epsilon\sum_{k\ge k_\epsilon}\mu_{k+1}$.
If $\omega$ were in the set $\Omega$, this would lead to a contradiction since the right-hand side diverges, due to the assumption $\sum_{k=0}^\infty\mu_k=\infty$, whereas the left-hand side is the remainder of a convergent series. Consequently, we have $\omega\notin\Omega$. Since $\bbP(\Omega^c)=0$, we thus obtain $\liminf_k\Psi(\xbs_k)=0$ a.s.
For every event in this a.s set $\Omega$ we can take a corresponding subsequence $(\xbs_{n_k})_{k\in\bbN}$ such that $\lim_{k\rightarrow\infty} \Psi(\xbs_{n_k})=0$.
Define now $\widehat\Psi(\xbs)=\sum_{i=1}^N \widehat\Psi_i(\xbs)$, with $\widehat\Psi_i(\xbs)=\yNi{F_i(\xbs)-\ybs_i}$.
We have $\liminf_k \widehat\Psi(\xbs_k)=0$ and $\lim_{k\rightarrow\infty} \widehat\Psi(\xbs_{n_k})=0$ (on the same subsequence).
Moreover, $\widehat\Psi(\xbs)^p \leq pN^p \Psi(\xbs)$.
The following argumentation is understood pointwise on the a.s. set $\Omega$ where $(\Delta_k)_{k\in\bbN}$ is convergent, $\sum_{k=0}^\infty \mu_{k+1}\Psi(\xbs_k)<\infty$, and $\liminf_k \widehat\Psi(\xbs_k)=0$.
Since $(\Delta_k)_{k\in\bbN}$ converges, it is bounded.
Note that $(\xbs_k)_{k\in\bbN}$ and $(\dmapX{p}(\xbs_k))_{k\in\bbN}$ are bounded by the coercivity of the Bregman distance, cf. Lemma \ref{lem:bregman_coercivity}.
By further passing to a subsequence, we can find a subsequence such that
$(\xN{\xbs_{n_k}})_{k\in\bbN}$ is convergent, $(\dmapX{p}(\xbs_{n_k}))_{k\in\bbN}$ is weakly convergent, and
\begin{align} \label{eqn:monotonic_subseq}
\lim_{k\rightarrow\infty} \widehat\Psi(\xbs_{n_k})=0\quad \text{and}\quad \widehat\Psi(\xbs_{n_k}) \leq \widehat\Psi(\xbs_{n}),\quad \forall n <n_k.
\end{align}
The latter can be obtained by setting $n_1=1$, and then recursively defining  $n_{k+1}=\min\{k> n_k: \widehat\Psi(\xbs_k)\leq \widehat\Psi(\xbs_{n_k})/2 \}$, for $k\in\bbN$. Any following subsequence satisfies the same property.
From \eqref{eqn:3_point_id}, for $k>\ell$ we have
\begin{align*}
    \bregman{\xbs_{n_\ell}}{\xbs_{n_k}}
    =&\frac{1}{p^\ast} \Big(\xN{\xbs_{n_\ell}}^p-\xN{\xbs_{n_k}}^p\Big) + \DP{\dmapX{p}(\xbs_{n_k})-\dmapX{p}(\xbs_{n_\ell})}{\xref}\\
    &+\DP{\dmapX{p}(\xbs_{n_k})-\dmapX{p}(\xbs_{n_\ell})}{\xbs_{n_k}-\xref}.
\end{align*}
Since the first two terms are Cauchy sequences, it suffices to treat the last term, denoted by ${\rm I}_{k,\ell}$.
Using telescopic sum, the iterate update rule and the Cauchy-Schwarz inequality, we have
\begin{align*}
    {\rm I}_{k,\ell}
    &=\sum_{n=n_\ell}^{n_k-1} \DP{\dmapX{p}(\xbs_{n+1})-\dmapX{p}(\xbs_{n})}{\xbs_{n_k}-\xref}\\
    &=\sum_{n=n_\ell}^{n_k-1} \mu_{n+1}\DP{F_{i_{n+1}}^{\prime}(\xbs_n)^\ast\svaldmap{p}(F_{i_{n+1}}(\xbs_n)-\ybs_{i_{n+1}})}{\xbs_{n_k}-\xref}\\
        &=\sum_{n=n_\ell}^{n_k-1} \mu_{n+1}\DP{\svaldmap{p}(F_{i_{n+1}}(\xbs_n)-\ybs_{i_{n+1}})}{F_{i_{n+1}}^\prime(\xbs_n)(\xbs_{n_k}-\xref)}\\
        &\leq\sum_{n=n_\ell}^{n_k-1} \mu_{n+1}\yNi{\svaldmap{p}(F_{i_{n+1}}(\xbs_n)-\ybs_{i_{n+1}})}\yNi{F_{i_{n+1}}^\prime(\xbs_n)(\xbs_{n_k}-\xref)}.
\end{align*}
Since $\widehat\Psi_i(\xbs)\leq \widehat\Psi(\xbs)$ holds for all $i\in[N]$, Using \eqref{eqn:TCC_1} and \eqref{eqn:monotonic_subseq} gives
\begin{align*}
    \yNi{F_{i_{n+1}}^\prime(\xbs_n)(\xbs_{n_k}-\xref)} &\leq \yNi{F_{i_{n+1}}^\prime(\xbs_n)(\xbs_{n}-\xref)}+\yNi{F_{i_{n+1}}^\prime(\xbs_n)(\xbs_{n_k}-\xbs_n)} \\
    &\leq (1+\gamma)\Big(\yNi{F_{i_{n+1}}(\xbs_{n})-\ybs_{i_{n+1}}}+\yNi{F_{i_{n+1}}(\xbs_{n_k})-F_{i_{n+1}}(\xbs_{n})}\Big)\\
    &\leq (1+\gamma)\Big(2\widehat\Psi_{i_{n+1}}(\xbs_{n})+\widehat\Psi_{i_{n+1}}(\xbs_{n_k})\Big)\leq 3(1+\gamma)\widehat\Psi(\xbs_{n}).
\end{align*}
Plugging this back into ${\rm I}_{k,\ell}$ yields
\begin{align*}
|{\rm I}_{k,\ell}|&\leq3(1+\gamma)\sum_{n=n_\ell}^{n_k-1} \mu_{n+1}\yNi{F_{i_{n+1}}(\xbs_n)-\ybs_{i_{n+1}}}^{p-1}\widehat\Psi(\xbs_{n})
    \leq3(1+\gamma)\sum_{n=n_\ell}^{n_k-1} \mu_{n+1}\widehat\Psi(\xbs_{n})^{p}.
\end{align*}
Since $\widehat\Psi(\xbs)^{p}\le pN^p \Psi(\xbs)$, the right-hand side of the above inequality converges to $0$ as $n_\ell\to\infty$.
Therefore, by \cite[Theorem 2.12(e)]{SLS_06}, it follows that $(\xbs_{n_k})_{k\in\bbN}$ is a Cauchy sequence, and thus it converges strongly to an $\widehat\xbs\in\CX_{\min}$.
Moreover, since $\CB_\nu(\xref)$ is closed we have $\widehat\xbs\in\CB_\nu(\xref)$.

The above argument showing the a.s. convergence of
$(\Delta_k)_{k\in\bbN}$ can be applied pointwise to any solution $\widehat\xbs\in\CX_{\min}$.
Namely, on the event where $(\xbs_{n_k})_{k\in\bbN}$ converges strongly to a solution $\widehat\xbs$, define $\widehat\Delta_k:=\bregman{\xbs_k}{\widehat\xbs}$. Similar to \eqref{eqn:bregman_descent}, we deduce
\begin{align*}
\widehat\Delta_{k+1}\leq\widehat\Delta_k -p\Big(1-\gamma -\Lmax^{p^\ast}\frac{G_{p^\ast}}{p^\ast}\mu_{k+1}^{p^\ast-1}\Big)\mu_{k+1}\Psi_{i_{k+1}}(\xbs_k).
\end{align*}
It follows that the (deterministic) sequence $(\widehat\Delta_k)_{k\in\bbN}$ converges to a $\widehat\Delta_\infty\geq0$.
The continuity of the Bregman distance in the first argument of Theorem \ref{thm:bregman_properties}(iv) gives
$\lim_{j\rightarrow\infty}\bregman{\xbs_{n_j}}{\widehat\xbs}=\bregman{\widehat\xbs}{\widehat\xbs}=0$,
and thus $\widehat\Delta_\infty=0$.
Moreover, by the $p$-convexity of $\CX$, cf. Theorem \ref{thm:bregman_properties}(iii), we have
$0\leq \xN{\xbs_k-\widehat\xbs}^p \leq \frac{p}{C_p} \widehat\Delta_k.$
From the squeeze theorem it follows that  $$\lim_{k\rightarrow\infty}\xN{\xbs_k-\widehat\xbs}=0.$$
Thus, for every event in an a.s. set $\Omega$, the sequence $(\xbs_k)_{k\in\bbN}$ strongly converges to some minimizing solution, that is
\begin{align*}
\bbP\Big(\lim_{k\rightarrow\infty} \inf_{\widetilde \xbs\in \CX_{\min}}\xN{\xbs_k-\widetilde \xbs}=0\Big)=1.
\end{align*}

We now turn to showing that the SGD iterates $(\xbs_k)_{k\in\bbN}$ converge a.s. to the $\xbs_0$-minimum-distance solution $\xref$.
By the construction of the scheme \eqref{eqn:sgd}, for all $k\in\bbN$, we have $\dmapX{p}(\xbs_{k+1})-\dmapX{p}(\xbs_k)\in\overline{\mathrm{range}(F^\prime(\xbs)^\ast)}$.
Using the condition $\Null(F^\prime(\xref))\subset\Null(F^\prime(\xbs))$ gives $\dmapX{p}(\xbs_{k+1})-\dmapX{p}(\xbs_k)\in \Null(F^\prime(\xref))^\perp$.
By applying a telescopic sum, we conclude that $\dmapX{p}(\xbs_{k})-\dmapX{p}(\xbs_0)\in \Null(F^\prime(\xref))^\perp$ for all $k\in\bbN$.
Since $\Null(F^\prime(\xref))^\perp$ is closed and $\dmapX{p}$ is continuous, we have $\dmapX{p}(\widehat{\xbs})-\dmapX{p}(\xbs_0)\in \Null(F^\prime(\xref))^\perp$.
From \eqref{eqn:TCC_1}, we get $F^\prime(\xref)(\widehat\xbs-\xref)=0$, i.e., $\widehat\xbs-\xref\in\Null(F^\prime(\xref))$. Therefore,
\begin{equation*}
    \DP{\dmapX{p}(\widehat{\xbs})-\dmapX{p}(\xbs_0)}{\widehat\xbs-\xref}=0.
\end{equation*}
Since $\widehat{\xbs}\in\CB_\nu(\xref)\subset\goodset$, the above identity and Lemma \ref{lem:unique_mini_dis_solu} imply $\DP{\dmapX{p}(\widehat{\xbs})-\dmapX{p}(\xref)}{\widehat\xbs-\xref}\le0$.
Theorem \ref{thm:dmap_properties}(i) then gives $\DP{\dmapX{p}(\widehat{\xbs})-\dmapX{p}(\xref)}{\widehat\xbs-\xref}=0$ and thus $\widehat\xbs=\xref$, since $\CX$ is $p$-convex.
\end{proof}

\begin{theorem}\label{thm:app_cone_convergence_qeqp}
Let the conditions of Lemma \ref{lem:bregman_monotonicity_qeqp} hold and that problem \eqref{eqn:inv_F_i} admits a solution $\xbs^\ast$.
If $\sum_{k=1}^\infty\mu_k=\infty$, then the sequence $(\xbs_k)_{k\in\bbN}$ converges a.s. to a minimizer $\widehat\xbs$ of $\Psi(\xbs)$.
Moreover, if $\overline{\mathrm{range}(F^\prime(\xbs)^\ast)} \subseteq \overline{\mathrm{range}(F^\prime(\xref)^\ast)}$ for all $\xbs \in \CB_\nu(\xref)$ and $\dmapX{p}(\xbs_0)\in \overline{\mathrm{range}(F^\prime(\xref)^\ast)}$, then $(\xbs_k)_{k\in\bbN}$ converges a.s. to the $\xbs_0$-minimum-norm solution $\xref$: $\lim_{k\to\infty}\bregman{\xbs_k}{\xref}=0$ a.s.
\end{theorem}
\begin{proof}
The proof that the sequence $(\xbs_k)_{k\in\bbN}$ converges a.s. to a minimizer $\widehat\xbs$ is identical to that of Theorem \ref{thm:cone_convergence_qeqp}. It suffices to focus on the second part.
Due to $\overline{\mathrm{range}(F^\prime(\xbs)^\ast)} \subseteq \overline{\mathrm{range}(F^\prime(\xref)^\ast)}$ for all $\xbs \in \CB_\nu(\xref)$ and $\dmapX{p}(\xbs_0)\in \overline{\mathrm{range}(F^\prime(\xref)^\ast)}$, it follows from \eqref{eqn:sgd} that $\dmapX{p}(\xbs_k)$ remains in $\overline{\mathrm{range}(F^\prime(\xref)^\ast)}$ for all $k\ge1$. The continuity of $\dmapX{p}$ then yields $\dmapX{p}(\widehat\xbs)\in\overline{\mathrm{range}(F^\prime(\xref)^\ast)}$, where $\widehat\xbs$ is some  minimizer of $\Psi(\xbs)$. By \eqref{eqn:TCC_1}, we have $F^\prime(\xref)(\widehat\xbs-\xref)=0$. Thus, from \cite[Proposition 7.1]{TBHK_12}, it follows that $\widehat\xbs=\xref$.
\end{proof}

\section{The proof of Theorem \ref{thm:L1_cone_convergence_qeqp}}\label{app:L1_cone_convergence_qeqp}
\begin{proof}
By Lemma \ref{lem:bregman_monotonicity_qeqp}, $(\bregman{\xbs_k}{\xref})_{k\in\bbN}$ is bounded, and is thus uniformly integrable, and by Theorem \ref{thm:cone_convergence_qeqp}, it converges a.s. to $0$. Then, by Vitali's convergence theorem \cite[Theorem 4.5.4]{B07}, $(\Delta_k)_{k\in\bbN}$ converges to $0$ also in expectation.
Now by the $p$-convexity of $\CX$ and the monotonicity of expectation, we have
\begin{align*}
    0\leq\frac{C_p}{p}\lim_{k\rightarrow\infty}\bbE[\xN{\xbs_k-\xref}^p] \leq \lim_{k\rightarrow\infty}\bbE[\bregman{\xbs_k}{\xref}]=0.
\end{align*}
By the continuity of the power function and the Lyapunov inequality for $1\leq r {\leq p}$, we have
\begin{align*}
0\leq\lim_{k\rightarrow\infty}\bbE[\xN{\xbs_k-\xref}^r]\leq\lim_{k\rightarrow\infty}(\bbE[\xN{\xbs_k-\xref}^p])^{r/p}=0.
\end{align*}
To prove the last claim we use uniform smoothness of $\CX$ and \cite[Lemma 5.16]{TBHK_12}, to deduce
\begin{align*}
    \xsN{\dmapX{p}(\xbs_k)-\dmapX{p}(\xref)}^{p^\ast}\leq C \max\{1, \xN{\xbs_k}, \xN{\xref}\}^{p}\, {\overline{\rho}_\CX}(\xN{\xbs_k-\xref})^{p^\ast},
\end{align*}
where $\overline{\rho}_\CX(\tau)=\rho_\CX(\tau)/\tau$ is a modulus of smoothness function such that $\overline{\rho}(\tau)\leq1$ and $\lim_{\tau\rightarrow0}\overline{\rho}(\tau)=0$, cf. Definition \ref{defn:smoothness_and_convexity}.
By Lemmas \ref{lem:bregman_coercivity} and \ref{lem:bregman_monotonicity_qeqp},
$(\xN{\xbs_k}^p)_{k\in\bbN}$ is {(uniformly) bounded, and thus the sequence $(\xsN{\dmapX{p}(\xbs_k)-\dmapX{p}(\xref)}^{p^\ast})_{k\in\bbN}$ is bounded and} uniformly integrable.
Since $\lim_{k\rightarrow\infty}\bbE[\xN{\xbs_k-\xref}]=0$, it follows that $\xN{\xbs_k-\xref}$ converges to $0$ in probability, and thus by the continuous mapping theorem  ${\overline{\rho}_\CX}(\xN{\xbs_k-\xref})^{p^\ast}$ also converges to $0$ in probability.
By Vitaly's theorem, and uniform integrabity of the sequence $(\xsN{\dmapX{p}(\xbs_k)-\dmapX{p}(\xref)}^{p^\ast})_{k\in\bbN}$,  it converges to $0$ in measure, and the claim follows.
\end{proof}

\bibliographystyle{abbrv}
\bibliography{biblio}

\end{document}